\documentclass{amsart}
\usepackage{amsopn}
\usepackage{amssymb, amscd}

\newcommand{\nc}{\newcommand}

\nc{\vg}{\mathfrak{v} } \nc{\wg}{\mathfrak{w} }
\nc{\zg}{\mathfrak{z} } \nc{\ngo}{\mathfrak{n} }
\nc{\ngoq}{\mathfrak{n}^{\QQ} } \nc{\ngoz}{\mathfrak{n}^{\ZZ} }
\nc{\ggoq}{\mathfrak{g}^{\QQ} } \nc{\kg}{\mathfrak{k} }
\nc{\mg}{\mathfrak{m} } \nc{\bg}{\mathfrak{b} }
\nc{\ggo}{\mathfrak{g} } \nc{\ggob}{\overline{\mathfrak{g}} }
\nc{\sog}{\mathfrak{so} } \nc{\sug}{\mathfrak{su} }
\nc{\spg}{\mathfrak{sp}} \nc{\slg}{\mathfrak{sl} }
\nc{\glg}{\mathfrak{gl} } \nc{\cg}{\mathfrak{c} }
\nc{\rg}{\mathfrak{r} } \nc{\hg}{\mathfrak{h} }
\nc{\tg}{\mathfrak{t} } \nc{\ug}{\mathfrak{u} }
\nc{\dg}{\mathfrak{d} } \nc{\ag}{\mathfrak{a} }
\nc{\pg}{\mathfrak{p} } \nc{\sg}{\mathfrak{s} }
\nc{\lgo}{\mathfrak{l} } \nc{\fg}{\mathfrak{f} }

\nc{\pca}{\mathcal{P}} \nc{\nca}{\mathcal{N}} \nc{\lca}{\mathcal{L}}
\nc{\oca}{\mathcal{O}} \nc{\mca}{\mathcal{M}} \nc{\tca}{\mathcal{T}}
\nc{\aca}{\mathcal{A}}

\nc{\vp}{\varphi} \nc{\ddt}{\frac{{\rm d}}{{\rm d}t}}
\nc{\im}{\mathtt{i}}

\renewcommand{\Im}{{\rm Im}}
\renewcommand{\l}{\lambda}
\nc{\ala}{Anosov Lie algebra} \nc{\alas}{Anosov Lie algebras}

\nc{\SO}{{\mathrm SO}} \nc{\Spe}{{\mathrm Sp}} \nc{\Sl}{{\mathrm
SL}} \nc{\SU}{{\mathrm SU}} \nc{\Or}{{\mathrm O}} \nc{\U}{{\mathrm
U}} \nc{\Gl}{{\mathrm GL}} \nc{\Se}{{\mathrm S}} \nc{\Cl}{{\mathrm
Cl}} \nc{\Spin}{{\mathrm Spin}} \nc{\Pin}{{\mathrm Pin}}

\nc{\RR}{{\Bbb R}} \nc{\HH}{{\Bbb H}} \nc{\CC}{{\Bbb C}}
\nc{\ZZ}{{\Bbb Z}} \nc{\FF}{{\Bbb F}} \nc{\NN}{{\Bbb N}}
\nc{\QQ}{{\Bbb Q}} \nc{\PP}{{\Bbb P}}

\nc{\vs}{\vspace{.5cm}}
 \nc{\ip}{\langle\cdot,\cdot\rangle}
\nc{\la}{\langle} \nc{\ra}{\rangle} \nc{\unm}{\frac{1}{2}}
\nc{\unc}{\frac{1}{4}} \nc{\und}{\frac{1}{16}} \nc{\f}{\frac}

\nc{\no}{\vs\noindent} \nc{\lam}{\Lambda^2\ggo^*\otimes\ggo}
\nc{\tang}{{\rm T}} \nc{\dif}{{\rm d}} \nc{\preq}{\simeq_K}
\nc{\lb}{[\,,\,]}

\nc{\He}{\operatorname{Hess}} \nc{\ad}{\operatorname{ad}}
\nc{\Ad}{\operatorname{Ad}} \nc{\rank}{\operatorname{rank}}
\nc{\Irr}{\operatorname{Irr}} \nc{\End}{\operatorname{End}}
\nc{\Aut}{\operatorname{Aut}} \nc{\Inn}{\operatorname{Inn}}
\nc{\Der}{\operatorname{Der}} \nc{\Ker}{\operatorname{Ker}}
\nc{\Iso}{\operatorname{I}} \nc{\Diff}{\operatorname{Diff}}
\nc{\Lie}{\operatorname{L}} \nc{\tr}{\operatorname{tr}}

\nc{\degr}{\operatorname{dgr}} \nc{\sen}{\operatorname{sen}}
\nc{\modu}{\operatorname{mod}} \nc{\Ric}{\operatorname{Ric}}
\nc{\sym}{\operatorname{sym}} \nc{\sca}{\operatorname{sc}}
\nc{\scalar}{{\sf s}} \nc{\grad}{\operatorname{grad}}
\nc{\ricci}{\operatorname{ric}} \nc{\Rin}{\operatorname{M}}
\nc{\Le}{\operatorname{L}}
\nc{\level}{\operatorname{level}} \nc{\rad}{\operatorname{r}}
\nc{\abel}{\operatorname{ab}} \nc{\Pf}{\operatorname{Pf}}

\newtheorem{theorem}{Theorem}[section]
\newtheorem{proposition}[theorem]{Proposition}
\newtheorem{corollary}[theorem]{Corollary}
\newtheorem{lemma}[theorem]{Lemma}
\newtheorem{definition}[theorem]{Definition}
\newtheorem{remark}[theorem]{Remark}

\newtheorem{example}[theorem]{Example}

\title[Anosov diffeomorphisms on nilmanifolds]{Anosov diffeomorphisms
on nilmanifolds up to dimension $8$}

\author{Jorge Lauret,\quad Cynthia E. Will}

\address{FaMAF and CIEM, Universidad Nacional de C\'ordoba, Haya de la Torre s/n, 5000 C\'ordoba, Argentina}
\email{lauret@mate.uncor.edu, cwill@mate.uncor.edu}

\thanks{
Supported by CONICET fellowships and grants from FONCyT and
Fundaci\'on Antorchas (Argentina).}

\begin{document}

\maketitle

\begin{abstract}
After more than thirty years, the only known examples of Anosov
diffeomorphisms are hyperbolic automorphisms of infranilmanifolds.
 It is also important to note that the existence of an Anosov
automorphism is a really strong condition on an infranilmanifold.
Any Anosov automorphism determines an automorphism of the rational
Lie algebra determined by the lattice, which is hyperbolic and
unimodular (and conversely ...). These two conditions together are
strong enough to make of such rational nilpotent Lie algebras
(called Anosov Lie algebras) very distinguished objects. In this
paper, we classify Anosov Lie algebras of dimension less or equal
than 8.

As a corollary, we obtain that if an infranilmanifold of dimension
$n\leq 8$ admits an Anosov diffeomorphism $f$ and it is not a torus
or a compact flat manifold (i.e. covered by a torus), then n=6 or 8
and the signature of $f$ necessarily equals $\{ 3,3\}$ or $\{
4,4\}$, respectively.  We had to study the set of all rational forms
up to isomorphism of many real Lie algebras, which is a subject on
its own and it is treated in a section completely independent of the
rest of the paper.
\end{abstract}

\section{Introduction}\label{intro}

A diffeomorphism $f$ of a compact differentiable manifold $M$ is
called {\it Anosov} if it has a global hyperbolic behavior, i.e. the
tangent bundle $\tang M$ admits a continuous invariant splitting
$\tang M=E^+\oplus E^-$ such that $\dif f$ expands $E^+$ and
contracts $E^-$ exponentially.  These diffeomorphisms define very
special dynamical systems and it is then a natural problem to
understand which are the manifolds supporting them (see \cite{Sml}).
After more than thirty years, the only known examples are hyperbolic
automorphisms of infranilmanifolds (called {\it Anosov
automorphisms}) and it is conjectured that any Anosov diffeomorphism
is topologically conjugate to one of these (see \cite{Mrg}).  The
conjecture is known to be true in many particular cases: J. Franks
\cite{Frn} and A. Manning \cite{Mnn} proved it for Anosov
diffeomorphisms on infranilmanifolds themselves; Y. Benoist and F.
Labourie \cite{BnsLbr} in the case the distributions $E^+,E^-$ are
differentiable and the Anosov diffeomorphism preserves an affine
connection (for instance a symplectic form); and J. Franks
\cite{Frn} when $\dim{E^+}=1$ (see also \cite{Grm} for expanding
maps).  Since Anosov automorphisms have many additional dynamical
properties (see \cite{Vrj}), a general resolution of the conjecture
would be of great relevance.

It is also important to note that the existence of an Anosov
automorphism is a really strong condition on an infranilmanifold.
An infranilmanifold is a quotient $N/\Gamma$, where $N$ is a
nilpotent Lie group and $\Gamma\subset K\ltimes N$ is a lattice
(i.e. a discrete cocompact subgroup) which is torsion-free and $K$
is a compact subgroup of $\Aut(N)$.  Among some other more technical
obstructions (see \cite{Mlf2} for further information), the first
natural obstruction for the infranilmanifold $N/\Gamma$ to admit an
Anosov automorphism is that the nilmanifold $N/(\Gamma\cap N)$,
which is a finite cover of $N/\Gamma$, has to do so.

In the case of a nilmanifold $N/\Gamma$ (i.e. when $\Gamma\subset
N$), any Anosov automorphism determines an automorphism $A$ of the
rational Lie algebra $\ngoq=\Gamma\otimes\QQ$, the Lie algebra of
the rational Mal'cev completion of $\Gamma$, which is {\it
hyperbolic} (i.e. $|\lambda|\ne 1$ for any eigenvalue $\lambda$ of
$A$) and {\it unimodular} (i.e. $[A]_{\beta}\in\Gl_n(\ZZ)$ for some
basis $\beta$ of $\ngoq$). Recall that $\ngoq$ is a rational form of
the Lie algebra $\ngo$ of $N$.  These two conditions together are
strong enough to make of such rational nilpotent Lie algebras
(called {\it Anosov Lie algebras}) very distinguished objects.  It
is proved in \cite{Ito} and \cite{Dkm} that if $\Gamma_1$ and
$\Gamma_2$ are commensurable (i.e.
$\Gamma_1\otimes\QQ\simeq\Gamma_2\otimes\QQ$) then $N/\Gamma_1$
admits an Anosov automorphism if and only if $N/\Gamma_2$ does. All
this suggests that the class of rational Anosov Lie algebras is the
key algebraic structure to study if one attempts to classify
infranilmanifolds admitting an Anosov diffeomorphism.

Finally, if one is interested in just those Lie groups which are
simply connected covers of such infranilmanifolds, then the objects
to be studied are real nilpotent Lie algebras $\ngo$ supporting a
hyperbolic automorphism $A$ such that $[A]_{\beta}\in\Gl_n(\ZZ)$ for
some $\ZZ$-basis $\beta$ of $\ngo$ (i.e. with integer structure
constants).  Such Lie algebras will also be called {\it Anosov}. We
note that a real Lie algebra is Anosov if and only if it has an
Anosov rational form.

The following would be then a natural program to classify all the
infranilmanifolds up to homeomorphism of a given dimension $n$ which
admits an Anosov diffeomorphism:
\begin{itemize}
\item[(i)] Find all $n$-dimensional Anosov Lie algebras over $\RR$.

\item[(ii)] For each real Lie algebra $\ngo$ obtained in (i), determine which
rational forms of $\ngo$  are Anosov.

\item[(iii)] For each rational Lie algebra $\ngoq$ from (ii), classify up to
isomorphism all the lattices $\Gamma$ in $N$, the nilpotent Lie
group with Lie algebra $\ngoq\otimes\RR$, such that
$\Gamma\otimes\QQ=\ngoq$. In other words, classify up to isomorphism
all the lattices in the commensurability class corresponding to
$\ngoq$.

\item[(iv)] Given a nilmanifold $N/\Gamma$ from (iii), decide which of the finitely
many infranilmanifolds $N/\Lambda$ essentially covered by $N/\Gamma$
(i.e. $\Lambda\cap N\simeq\Gamma$) admits an Anosov automorphism,
that is, a hyperbolic automorphism $\vp$ of $N$ such that
$\vp(\Lambda)=\Lambda$ (see \cite{Mlf2}).
\end{itemize}

Parts (i) and (ii) have been solved for dimension $n\leq 6$ in
\cite{CssKnnScv} and \cite{Mlf}, yielding only two algebras over
$\RR$: $\hg_3\oplus\hg_3$ and $\fg_3$ (see Table \ref{notation}).
There are some other families of real Anosov Lie algebras in the
literature (see Remark \ref{exaAnosov}).  Besides these examples in
somewhat sporadic dimensions, there is a construction in
\cite{anosov} proving that $\ngo\oplus\ngo$ is Anosov for any real
graded nilpotent Lie algebra $\ngo$ which admits at least one
rational form $\ngoq$.  We note that the existing Anosov rational
form is not necessarily $\ngoq\oplus\ngoq$.  Since for instance any
$2$-step nilpotent Lie algebra is graded, this construction shows
that part (i) of the program above is already a wild problem for $n$
large.  Furthermore, by using the classification of nilpotent Lie
algebras in low dimensions (see \cite{Mgn,Sly}), we can assert that
there are at least $18$ real Anosov Lie algebras of dimension $10$,
$68$ in dimension $12$ and more than 100, together with some curves
in dimension $14$.

In view of this fact, the aim of this paper is to approach the
classification in small dimensions.  We have classified up to
isomorphism real and rational Lie algebras of dimension $\leq 8$
which are Anosov.  In other words, we have solved parts (i) and (ii)
of the program for $n=7$ and $n=8$.  We refer to Tables
\ref{notation} and \ref{Anosovtable} for a quick look at the results
obtained. Without an abelian factor, there are only three
$8$-dimensional real Lie algebras which are Anosov and none in
dimension $7$.  This is a really small list, bearing in mind that
there exist several one and two-parameters families and hundreds of
isolated examples of $7$ and $8$-dimensional nilpotent Lie algebras,
and there is not even a full classification in dimension $8$.

One of the corollaries which might be interesting from a dynamical
point of view is that if an infranilmanifold of dimension $n\leq 8$
admits an Anosov diffeomorphism $f$ and it is not a torus or a
compact flat manifold (i.e. covered by a torus), then $n=6$ or $8$
and the signature of $f$, defined by $\{\dim{E^+},\dim{E^-}\}$,
necessarily equals $\{ 3,3\}$ or $\{ 4,4\}$, respectively.

We now give an idea of the structure of the proof.  The {\it type}
of a nilpotent Lie algebra $\ngo$ is the $r$-tuple $(n_1,...,n_r)$ ,
where $n_i=\dim{C^{i-1}(\ngo)/C^i(\ngo)}$ and $C^i(\ngo)$ is the
central descending series. By using that any Anosov Lie algebra
admits an Anosov automorphism $A$ which is semisimple and some
 elementary properties of lattices, one sees that only a few types
 are allowed in each dimension $7$ and $8$.  We
then study these types case by case in Section \ref{cla} and exploit
that the eigenvalues of $A$ are algebraic integers (even units). For
each of the types we get only one or two real Lie algebras
(sometimes no one at all) which are candidates to be Anosov, and
some of them are excluded by using a criterion given in terms of a
homogeneous polynomial (called the {\it Pfaffian form}) associated
to each $2$-step nilpotent Lie algebra.

We previously study the set of all rational forms up to isomorphism
for each of the real Lie algebras obtained in the classification
over $\RR$. This is a subject on its own, and it is treated in
Section \ref{rational}, a part which is completely independent of
the rest of the paper.  The results obtained there (see Table
\ref{ratforms}) allows us to classify Anosov Lie algebras over $\QQ$
in Section \ref{cla2}, and here we also use a criterion on the
Pfaffian form to discard some of them, which has in this case
integer coefficients and hence some topics from number theory as the
Pell equation and square free numbers appear.  Such criterions and
most of the known tools to deal with Anosov automorphisms are given
in Section \ref{andiff} (see also \cite{Dn} for an approach via
representation theory and arithmetic groups), as well as a
generalization of the construction in \cite{anosov} suggested by F.
Grunewald, proving that $\ngo\oplus...\oplus\ngo$ ($s$ times, $s\geq
2$) is Anosov for any graded nilpotent Lie algebra over $\RR$ having
a rational form.

\vs \noindent {\it Acknowledgements.}  We wish to thank M. Mainkar
and S.G. Dani for very helpful comments on a first version of this
paper.

\section{Rational forms of nilpotent Lie algebras}\label{rational}

Since the classification of all nilmanifolds admitting an Anosov
diffeomorphism reduces to the determination of a special class of
nilpotent Lie algebras over $\QQ$, we now start the study of
rational forms of real nilpotent Lie algebras.  Let $\ngo$ be a
nilpotent Lie algebra over $\RR$ of dimension $n$.

\begin{definition}\label{ratform}
{\rm A {\it rational form} of $\ngo$ is an $n$-dimensional rational
subspace $\ngoq$ of $\ngo$ such that
$$
[X,Y]\in\ngoq, \qquad\forall\; X,Y\in\ngoq.
$$
Two rational forms $\ngoq_1$, $\ngoq_2$ of $\ngo$ are said to be
{\it isomorphic} if there exists $A\in\Aut(\ngo)$ such that
$A\ngoq_1=\ngoq_2$, or equivalently, if they are isomorphic as Lie
algebras over $\QQ$ (recall that $\ngoq\otimes\RR=\ngo)$.
Analogously, by considering $\RR$ and $\CC$ (resp. $\QQ$ and $\CC$)
instead of $\QQ$ and $\RR$ one defines a {\it real form} (resp. a
{\it rational form}) of a complex Lie algebra.  }
\end{definition}

Not every real nilpotent Lie algebra admits a rational form.  By a
result due to Malcev, the existence of a rational form of $\ngo$ is
equivalent to the corresponding Lie group $N$ admits a {\it
lattice}, i.e. a cocompact discrete subgroup (see \cite{Rgn}).
Another difference with the semisimple case is that sometimes $\ngo$
has only one rational form up to isomorphism.  The problem of
finding all isomorphism classes of rational forms for a given real
nilpotent Lie algebra is a very difficult one, even in the low
dimensional or two-step cases. Very little is known about this
challenge problem in the literature (see \cite[Section 5]{Ebr}).
When $\ngo$ is two-step nilpotent and has $2$-dimensional center, F.
Grunewald and D. Segal \cite{GrnSgl1,GrnSgl2} gave an answer in
terms of isomorphism classes of binary forms, which will be
explained below.  In \cite{Smn} it is proved that
$\hg_{2k+1}\oplus\RR^m$ has only one rational form up to isomorphism
for all $k,m$, and that certain real Lie algebras of the form
$\ggo\oplus\ggo$ have infinitely many ones.  In this section, we
will find all the rational forms up to isomorphism of four real
nilpotent Lie algebras of dimension $8$ (see Table \ref{ratforms}).
This information will be useful in the classification of
$8$-dimensional Anosov Lie algebras.

\vs

Let $\ngo$ be a Lie algebra over the field $K$, which is assumed
from now on to be of characteristic zero.  We are mainly interested
in the cases $K=\CC,\RR,\QQ$. Fix a positive definite symmetric
$K$-bilinear form $\ip$ on $\ngo$ (i.e. an inner product).  For each
$Z\in\ngo$ consider the $K$-linear transformation
$J_Z:\ngo\mapsto\ngo$ defined by
\begin{equation}\label{jota}
\la J_ZX,Y\ra=\la [X,Y],Z\ra, \qquad\forall\; X,Y\in\ngo.
\end{equation}
Recall that $J_Z$ is skew symmetric with respect to $\ip$ and the
map $J:\ngo\mapsto\sog(n,K)$ is $K$-linear, where $n$ is the
dimension of $\ngo$. Equivalently, we may define these maps by
fixing a basis $\beta=\{ X_1,...,X_n\}$ of $\ngo$ rather than an
inner product in the following way: $J_Z$ is the $K$-linear
transformation whose matrix in terms of $\beta$ is
$$
\left(\sum_{k=1}^{n}c_{ij}^kx_k\right),
$$
where $[X_i,X_j]=\displaystyle{\sum_{k=1}^{n}}c_{ij}^kX_k$ and
$Z=\displaystyle{\sum_{k=1}^{n}}x_{k}X_k$. This coincides with the
first definition if one sets $\la X_i,X_j\ra=\delta_{ij}$.

If $\ngo$ and $\ngo'$ are two Lie algebras over $K$ and $\{ J_Z\}$,
$\{ J_Z'\}$ are the corresponding maps relative to the inner
products $\ip$ and $\ip'$ respectively, then it is not hard to see
that a linear map $A:\ngo\mapsto\ngo'$ is a Lie algebra isomorphism
if and only if
\begin{equation}\label{iso}
A^tJ_Z'A=J_{A^tZ}, \qquad\forall\; Z\in\ngo,
\end{equation}
where $A^t:\ngo'\mapsto\ngo$ is given by $\la A^tX,Y\ra=\la X,AY\ra$
for all $X\in\ngo'$, $Y\in\ngo$.

\begin{definition}\label{type}
{\rm Consider the central descendent series defined by
$C^0(\ngo)=\ngo$, $C^i(\ngo)=[\ngo,C^{i-1}(\ngo)]$.  When
$C^r(\ngo)=0$ and $C^{r-1}(\ngo)\ne 0$, $\ngo$ is said to be
$r$-step nilpotent, and we denote by $(n_1,...,n_r)$ the {\it type}
of $\ngo$, where
$$
n_i=\dim{C^{i-1}(\ngo)/C^i(\ngo)}.
$$
We also take a decomposition $\ngo=\ngo_1\oplus...\oplus\ngo_r$, a
direct sum of vector spaces, such that
$C^i(\ngo)=\ngo_{i+1}\oplus...\oplus\ngo_r$ for all $i$.}
\end{definition}

Assume now that $\ngo$ is $2$-step nilpotent, or equivalently of
type $(n_1,n_2)$. We will always have fixed orthonormal basis $\{
X_i\}$ and $\{ Z_j\}$ of $\ngo_1$ and $\ngo_2$, respectively.
Consider any direct sum decomposition of the form
$\ngo=V\oplus[\ngo,\ngo]$, that is, $\ngo_1=V$. If the inner product
satisfies $\la V,[\ngo,\ngo]\ra=0$ then $V$ is $J_Z$-invariant for
any $Z$ and $J_Z=0$ if and only if $Z\in V$. We define
$f:[\ngo,\ngo]\mapsto K$ by
$$
f(Z)=\Pf(J_Z|_V), \qquad Z\in [\ngo,\ngo],
$$
where $\Pf:\sog(V,K)\mapsto K$ is the {\it Pfaffian}, that is, the
only polynomial function satisfying $\Pf(B)^2=\det{B}$ for all
$B\in\sog(V,K)$ and $\Pf(J)=1$ for
$$
J= \left[\begin{array}{cc}
0&I\\
-I&0
\end{array}\right].
$$
Roughly speaking, $f(Z)=\left(\det{J_Z|_V}\right)^{\unm}$, and so we
need $\dim{V}$ to be even in order to get $f\ne 0$. For any
$A\in\glg(V,K)$, $B\in\sog(V,K)$ we have that
$\Pf(ABA^t)=(\det{A})\Pf(B)$.

\begin{definition}\label{pform}
{\rm We call $f$ the {\it Pfaffian form} of $\ngo$.}
\end{definition}

If $\dim{V}=2m$ and $\dim{[\ngo,\ngo]}=k$ then $f=f(x_1,...,x_k)$ is
a homogeneous polynomial of degree $m$ in $k$ variables with
coefficients in $K$, where $Z=\sum_{i=1}^{k}x_iZ_i$ and $\{
Z_1,...,Z_k\}$ is a fixed basis of $[\ngo,\ngo]$. $f$ is also called
a form of degree $m$, when $k=2$ or $3$ one uses the words binary or
ternary and for $m=2,3$ and $4$, quadratic, cubic and cuartic,
respectively.

Let $P_{k,m}(K)$ denote the set of all homogeneous polynomials of
degree $m$ in $k$ variables with coefficients in $K$.  The group
$\Gl_k(K)$ acts naturally on $P_{k,m}(K)$ by
$$
(A.f)(x_1,...,x_k)=f(A^{-1}(x_1,...,x_k)),
$$
that is, by linear substitution of variables, and thus the action
determines the usual equivalence relation between forms, denoted by
$f\simeq g$.  In the present paper, we need to consider the
following wider equivalence relation.

\begin{definition}\label{equiv} {\rm For $f,g\in P_{k,m}(K)$, we say that $f$ is {\it projectively equivalent}
to $g$, and denote it by $f\preq g$, if there exists $A\in\Gl_k(K)$
and $c\in K^*$ such that}
$$
f(x_1,...,x_k)=cg(A(x_1,...,x_k)).
$$
\end{definition}

In other words, we are interested in projective equivalence classes
of forms.

\begin{proposition}\label{isoforms}
Let $\ngo,\ngo'$ be two-step nilpotent Lie algebras over the field
$K$.  If $\ngo$ and $\ngo'$ are isomorphic then $f\preq f'$, where
$f$ and $f'$ are the Pfaffian forms of $\ngo$ and $\ngo'$,
respectively.
\end{proposition}

\begin{proof}
Since $\ngo$ and $\ngo'$ are isomorphic we can assume that
$\ngo=\ngo'$ and $[\ngo,\ngo]=[\ngo',\ngo']$ as vector spaces, and
then the decomposition $\ngo=V\oplus[\ngo,\ngo]$  is valid for both
Lie brackets $\lb$ and $\lb'$.  Any isomorphism satisfies
$A[\ngo,\ngo]=[\ngo',\ngo']'$, and it is easy to see that there is
always an isomorphism $A$ between them satisfying $AV=V$. It follows
from (\ref{iso}) that
$$
A^tJ_Z'A=J_{A^tZ}, \qquad \forall\; Z\in [\ngo,\ngo],
$$
and since the subspaces $V$ and $[\ngo,\ngo]$ are preserved by $A$
and $A^t$ we have that
$$
f'(Z)=cf(A_2^tZ),
$$
where $A_2=A|_{[\ngo,\ngo]}$ and $c^{-1}=\det{A|_V}$.  This shows
that $f\preq f'$.
\end{proof}

The above proposition says that the (projective) equivalence class
of the form $f(x_1,...,x_k)$ is an isomorphism invariant of the Lie
algebra $\ngo$.  This invariant was introduced by Scheuneman in
\cite{Sch}.

What is known about the classification of forms?.  Unfortunately,
much less than one could naively expect.  The case $K=\CC$ is as
usual the most developed one, and there the understanding of the
ring of invariant polynomials $\CC[P_{k,m}]^{\Sl_k(\CC)}$ is
crucial.  A set of generators and their relations for such a ring is
known only for small values of $k$ and $m$, for instance for $k=2$
and $m\leq 8$, or $k=3$ and $m\leq 3$.  The following well known
result will help us to distinguish between projective equivalence
classes of forms, and in view of Proposition \ref{isoforms}, to
recognize non-isomorphic two-step nilpotent Lie algebras.

\begin{proposition}\label{hessian}
If $f,g\in P_{k,m}(K)$ satisfy
$$
f(x_1,...,x_k)=cg(A(x_1,...,x_k))
$$
for some $A\in\Gl_k(K)$ and $c\in K^*$, then
$$
Hf(x_1,...,x_k)=c^k(\det{A})^2Hg(A(x_1,...,x_k)),
$$
where the Hessian $Hf$ of the form $f$ is defined by
$$
Hf(x_1,...,x_k)=\det{\left[\frac{\partial^2f}{\partial x_i\partial
x_j}\right]}.
$$
\end{proposition}

\vs

Let $\ngoq$ be a rational nilpotent Lie algebra of type $(4,2)$. If
$\ngoq=\ngo_1\oplus\ngo_2$ is the decomposition such that
$\dim{\ngo_1}=4$, $\dim{\ngo_2}=2$ and $[\ngoq,\ngoq]=\ngo_2$, then
we consider the Pfaffian form $f$ of $\ngo$.  Thus $f$ is a binary
quadratic form, say $f(x,y)=ax^2+bxy+cy^2$, with $a,b,c\in\QQ$. The
strong result proved in \cite{GrnSgl1} is that the converse of
Proposition \ref{isoforms} is valid in this case, that is, there is
a one-to-one correspondence between isomorphism classes of
non-degenerate (i.e. with center equal to $\ngo_2$) rational Lie
algebras of type $(4,2)$ and projective equivalence classes of
binary quadratic forms with coefficients in $\QQ$.  It is easy to
see that such classes can be parametrized by
$$
\{ f_k(x,y)=x^2-ky^2:k\;\;\mbox{is a square free integer number}\}.
$$
Recall that an integer number is said to be {\it square free} if
$p^2\nmid k$ for any prime $p$, and the set of all square free
numbers parametrizes the equivalence classes of the relation in
$\QQ$ defined by $r\equiv s$ if and only if $r=q^2s$ for some
$q\in\QQ^*$.  We are considering $k=0$ a square free number too. If
$f_k\preq f_{k'}$ then it follows from Proposition \ref{hessian}
that $-4k=-4q^2k'$ for some $q\in\QQ^*$, which implies that $k=k'$
in the case $k$ and $k'$ are square free.

It is not hard to prove that the Pfaffian form of the Lie algebra
$\ngoq_k=\ngo_1\oplus\ngo_2$ defined by
\begin{equation}\label{crath3h3}
[X_1,X_3]=Z_1, \quad [X_1,X_4]=Z_2, \quad [X_2,X_3]=kZ_2, \quad
[X_2,X_4]=Z_1
\end{equation}
is $f_k$.  For $K=\RR$, these Lie algebras can be distinguished only
by the sign of the discriminant of $f_k$, which says that there are
only three real Lie algebras of type $(4,2)$, namely, those of the
form $\ngoq_k\otimes\RR$ with $k>0$, $k=0$ and $k<0$. We have that
$\ngoq_1\otimes\RR\simeq\hg_3\oplus\hg_3$, where $\hg_3$ denotes the
$3$-dimensional Heisenberg Lie algebra and $\ngoq_{-1}\otimes\RR$ is
an H-type Lie algebra.  Analogously, there are only two
complexifications $\ngoq_k\otimes\CC$, those with $k\ne 0$ and
$k=0.$

\begin{proposition}\label{rath3h3}
The set of isomorphism classes of rational forms of the Lie algebra
$\hg_3\oplus\hg_3$ is parametrized by
$$
\{ \ngoq_k:k\,\mbox{is a square free natural number}\}.
$$
\end{proposition}

\begin{proof}
The Lie bracket of $\hg_3\oplus\hg_3$ is
$$
[X_1,X_2]=Z_1, \qquad [X_3,X_4]=Z_2,
$$
and one can easily check that the rational subspace generated by the
set
$$
\left\{ X_1+X_3, \; \sqrt{k}(X_1-X_3), \; \sqrt{k}(X_2+X_4), \;
X_2-X_4,\; \sqrt{k}(Z_1+Z_2), \; Z_1-Z_2\right\},
$$
is a rational subalgebra of $\hg_3\oplus\hg_3$ isomorphic to
$\ngoq_k$.
\end{proof}

We now describe the results in \cite{GrnSgl2} for the general case
(see also \cite{Ggr}).  Consider $\ngo=\ngo_1\oplus\ngo_2$ a vector
space over $K$ such that $\ngo_1$ and $\ngo_2$ are subspaces of
dimension $n$ and $2$ respectively. Every $2$-step nilpotent Lie
algebra of dimension $n+2$ and $2$-dimensional center can be
represented by a bilinear form $\mu:\ngo_1\times\ngo_1\mapsto\ngo_2$
which is non-degenerate in the following way: for any nonzero
$X\in\ngo_1$ there is a $Y\in\ngo_1$ such that $\mu(X,Y)\ne 0$.  If
we fix basis $\{ X_1,...,X_n\}$ and $\{ Z_1,Z_2\}$ of $\ngo_1$ and
$\ngo_2$ respectively, then each $\mu$ has an associated Pfaffian
binary form $f_{\mu}$ defined by
$$
f_{\mu}(x,y)=\Pf(J^{\mu}_{xZ_1+yZ_2})
$$
(see Definition \ref{pform}).  A {\it central decomposition} of
$\mu$ is given by a decomposition of $\ngo_1$ in a direct sum of
subspaces $\ngo_1=V_1\oplus...\oplus V_r$ such that $\mu(V_i,V_j)=0$
for all $i\ne j$.  We say that $\mu$ is {\it indecomposable} when
the only possible central decomposition has $r=1$.  Every $\mu$ has
a central decomposition into indecomposables constituents and such a
decomposition is unique up to an automorphism of $\mu$; in
particular, the constituents $V_i\oplus\ngo_2$ are unique up to
isomorphism.

There is only one indecomposable $\mu$ for $n$ odd and it can be
defined by
$$
J^{\mu}_{xZ_1+yZ_2}=\left[\begin{array}{c|c} 0&\begin{array}{ccccc}
-x&-y&&&0\\
0&-x&-y&&\\
&&\ddots&\ddots&\\
0&&&-x&-y \end{array} \\
\hline
\begin{array}{cccc}
x&0&&0 \\
y&x&& \\
0&y&\ddots& \\
&&\ddots&x \\
0&&&y
\end{array} & 0
\end{array}\right].
$$
Recall that $f_{\mu}=0$ in this case.  When $n$ is even the
situation is much more abundant: two indecomposables $\mu$ and
$\lambda$ are isomorphic if and only if $f_{\mu}\preq f_{\lambda}$.
If $n=2m$ and $f_{\mu}(x,y)=x^m-a_1x^{m-1}y-...-a_my^m$, then
$$
J^{\mu}_{xZ_1+yZ_2}=\left[\begin{array}{cc} 0&-B^t\\
B&0\end{array}\right],
$$
where
$$
B=\left[\begin{array}{ccccc}
x&y&0&\cdots &0\\
0&x&y&&\vdots\\
\vdots&&\ddots&\ddots&0\\
0&\cdots&0&x&y\\
a_my&a_{m-1}y&\cdots&a_2y&a_1y+x
\end{array}\right].
$$
We note that here $f_{\mu}$ is always nonzero, and in order to get
$\mu$ indecomposable one needs the form $f_{\mu}$ to be primitive
(i.e. a power of an irreducible one).  For decomposable $\mu$ and
$\lambda$ with respective central decompositions
$\ngo_1=V_1\oplus...\oplus V_r$ and $\ngo_1=W_1\oplus...\oplus W_s$
into indecomposables constituents, we have that $\mu$ is isomorphic
to $\lambda$ if and only if $r=s$ and after a suitable reordering
one has that
\begin{itemize}
\item[(i)] for some $t\leq r$, $\dim{V_i}=\dim{W_i}$ for all $i=1,...,t$ and they
are all even numbers;

\item[(ii)] if $\mu_i=\mu|_{V_i\times V_i}$, $\lambda_i=\lambda|_{W_i\times W_i}$
then there exist $A\in\Gl_2(K)$ and $c_1,...,c_t\in K^*$ such that
$$
f_{\mu_i}(x,y)=c_if_{\lambda_i}(A(x,y)) \qquad \forall\; i=1,...,t;
$$
\item[(iii)] $\dim{V_i}=\dim{W_i}$ is odd for all $i=t+1,...,r$.
\end{itemize}

Concerning our search for all rational forms up to isomorphism of a
given real nilpotent Lie algebra, these results say that the picture
in the $2$-step nilpotent with $2$-dimensional center case is as
follows.  Let $(\ngoq=\ngo_1\oplus\ngo_2, \mu)$ be one of such Lie
algebras over $\QQ$, and consider the corresponding Pfaffian form
$f_{\mu}\in P_{2,m}(\QQ)$.  The isomorphism classes of rational
forms of $\ngoq\otimes\RR$ are then parametrized by
$$
\Big((\RR^*\times\Gl_2(\RR)).f_{\mu}\cap
P_{2,m}(\QQ)\Big)/(\QQ^*\times\Gl_2(\QQ)).
$$
In other words, the rational points of the orbit
$(\RR^*\times\Gl_2(\RR)).f_{\mu}$ ($f_{\mu}$  viewed as an element
of $P_{2,m}(\RR)$) is a $(\QQ^*\times\Gl_2(\QQ))$-invariant set and
we have to consider the orbit space of this action.  Such a
description shows the high difficulty of the problem.  Recall that
we have to consider the action of $\RR^*\times\Gl_2(\RR)$ instead of
just that of $\Gl_2(\RR)$ only when $m$ is even.

In what follows, we will study rational forms of four
$8$-dimensional nilpotent Lie algebras.  We refer to Table
\ref{ratforms} for a summary of the results obtained.

Let $\ggo$ be the $8$-dimensional $2$-step nilpotent Lie algebra
defined by
\begin{equation}\label{62alg}
[X_1,X_2]=Z_1, \quad [X_1,X_3]=Z_2,\quad [X_4,X_5]=Z_1,\quad
[X_4,X_6]=Z_2.
\end{equation}
It is easy to see that its Pfaffian form $f$ is zero.  Let $\ggoq$
be a rational form of $\ggo$, for which we can assume that
$\ggoq=\la X_1,...,X_6\ra_{\QQ}\oplus\la Z_1,Z_2\ra_{\QQ}$.  Since
the Pfaffian form $g$ of $\ggoq$ satisfies $g\simeq_\RR f=0$ we
obtain that $g=0$.  By using the results described above one deduces
that $\ggoq$ can not be indecomposable, and so $\la
X_1,...,X_6\ra_{\QQ}=V_1\oplus...\oplus V_r$ with $[V_i,V_j]=0$ for
all $i\ne j$. Now, $\la
X_1,...,X_6\ra_{\RR}=V_1\otimes\RR\oplus...\oplus V_r\otimes\RR$ is
also a central decomposition for $\ggo$, proving that $r=2$ and
$\dim{V_1}=\dim{V_2}=3$ by the uniqueness of such a decomposition.
By applying again the results described above, now to the odd
situation, we get the following

{\small
\begin{table}
$$
\begin{array}{ccc}
\hline\hline
&& \\
Notation & Type & Lie \quad brackets\\ \\
\hline
&&\\
\hg_{2k+1}&(2k,1) & [X_1,X_2]=Z_1,..., [X_{2k-1},X_{2k}]=Z_1 \quad   \\ \\

\fg_3&(3,3) & [X_1,X_2]=Z_1,\,[X_1,X_3]=Z_2,\,[X_2,X_3]=Z_3  \\ \\

\ggo&(6,2) & [X_1,X_2]=Z_1,\,[X_1,X_3]=Z_2,\,[X_4,X_5]=Z_1,\,[X_4,X_6]=Z_2  \\ \\

\hg &(4,4) & [X_1,X_3]=Z_1,\,[X_1,X_4]=Z_2,\,[X_2,X_3]=Z_3,\,[X_2,X_4]=Z_4   \\ \\

\lgo_4 & (2,1,1) & [X_1,X_2]=X_3,\,[X_1,X_3]=X_4 \\ \\

\hline \hline \\
\end{array}
$$
\caption{Notation for some real nilpotent Lie
algebras.}\label{notation}
\end{table}}

\begin{proposition}\label{rat62} The Lie algebra $\ggo$ of type $(6,2)$ given in {\rm (\ref{62alg})} has only one
rational form up to isomorphism, denoted by $\ggoq$.
\end{proposition}

\begin{remark}\label{frg}
{\rm Clearly, the same proof is valid if one need to find all real
forms of the complex Lie algebra $\ggo_{\CC}=\ggo\otimes\CC$.  Thus
$\ggo$ is the only real form of $\ggo_{\CC}$ up to isomorphism.}
\end{remark}

As another application of the correspondence with binary forms given
above, we now study rational forms of the real Lie algebra
$\hg_3\oplus\hg_5$ of type $(6,2)$.  It has central decomposition
$\ngo_1=V_1\oplus V_2\oplus V_3$ with $\dim{V_i}=2$ for all $i$ as a
real Lie algebra and its Pfaffian form is $f(x,y)=xy^2$.  Let
$\mu:\ngo_1\times\ngo_1\mapsto\ngo_2$ be a rational form of
$\hg_3\oplus\hg_5$ with Pfaffian form $f_{\mu}$.  If $\mu$ is
decomposable then $\ngo_1=W_1\oplus W_2$, $\dim{W_1}=2$,
$\dim{W_2}=4$; or $\ngo_1=W_1\oplus W_2\oplus W_3$, $\dim{W_i}=2$
for all $i$.  In any case, $f_{\mu_i}\simeq_{\QQ}x,y$ or $y^2$
proving that $\mu$ must be isomorphic to the canonical rational form
$$
\mu_0(X_1,X_2)=Z_1, \quad \mu_0(X_3,X_4)=Z_2,
\quad\mu_0(X_5,X_6)=Z_2,
$$
for which $f_{\mu_0}=f$.  We then assume that $\mu$ is
indecomposable.  We shall prove that there is only one
$\Gl_2(\QQ)$-orbit of rational points in $\Gl_2(\RR).f$, and so
$\mu$ will have to be isomorphic to $\mu_0$.  There exists
$A\in\Gl_2(\RR)$ such that $f_{\mu}=A^{-1}.f$, that is,
$$
f_{\mu}(x,y)=ac^2x^3+c(2ad+bc)x^2y+d(ad+2bc)xy^2+bd^2y^3, \qquad A=\left[\begin{smallmatrix} a&b\\
c&d\end{smallmatrix}\right].
$$
Since $\mu$ is rational we have that
$$
q:=ac^2, \qquad r:=c(2ad+bc), \qquad s:=d(ad+2bc), \qquad t:=bd^2
$$
are all in $\QQ$.  If $c=0$ then $q=r=0$ and $s=ad^2$, $t=bd^2$,
which implies that $s\ne 0$ and hence
$$
f_{\mu}=B^{-1}.f, \qquad \mbox{for}\quad B=\left[\begin{smallmatrix} s&t\\
0&1\end{smallmatrix}\right]\in\Gl_2(\QQ).
$$
If $c\ne 0$ then one can check by a straightforward computation that
$$
\frac{d}{c}=\frac{9qst+rs^2-6r^2t}{6qs^2-r^2s-9qrt} \in\QQ.
$$
There must be a simpler formula for $\frac{d}{c}$ in terms of
$q,r,s,t$, but unfortunately we were not able to find it.  By
putting $u:=\frac{d}{c}$ we have that
$$
f_{\mu}=B^{-1}.f, \qquad \mbox{for} \quad B=\left[\begin{smallmatrix} q&t/u^2 \\
1&u\end{smallmatrix}\right]\in\Gl_2(\QQ).
$$
Recall that $\det{B}=qu-\frac{t}{u^2}=c(ad-bc)=c\det{A}\ne 0$.  We
then obtain that in any case $f_{\mu}\simeq_{\QQ}f$ and so $\mu$ is
isomorphic to $\mu_0$.

\begin{proposition}\label{rath3h5} Up to isomorphism, the real Lie algebra $\hg_3\oplus\hg_5$ of type $(6,2)$ has only one rational form, which will be denoted by $(\hg_3\oplus\hg_5)^{\QQ}$.
\end{proposition}

\begin{remark}\label{frh3h5}
{\rm It is easy to check that the above proof is also valid if we
replace $\QQ$ and $\RR$ by $\RR$ and $\CC$, obtaining in this way
that the only real form of $(\hg_3\oplus\hg_5)_{\CC}$ is
$\hg_3\oplus\hg_5$.}
\end{remark}

We now describe a duality for $2$-step nilpotent Lie algebras over
any field of characteristic zero introduced by J. Scheuneman
\cite{Sch} (see also \cite{Ggr}), which assigns to each Lie algebra
of type $(n,k)$ another one of type $(n,\frac{n(n-1)}{2}-k)$. The
dual of a Lie algebra $\ngo=\ngo_1\oplus\ngo_2$ of type $(n,k)$ can
be defined as follows: consider the maps $\{
J_Z:Z\in\ngo_2\}\subset\sog(n)$ corresponding to a fixed inner
product $\ip$ on $\ngo$ (see (\ref{jota})). Let
$\tilde{\ngo}_2\subset\sog(n)$ be the orthogonal complement of the
$k$-dimensional subspace $\{ J_Z:Z\in\ngo_2\}$ in $\sog(n)$ relative
to the inner product $(A,B)=-\tr{AB}$.  Now, we define the $2$-step
nilpotent Lie algebra $\tilde{\ngo}=\ngo_1\oplus\tilde{\ngo}_2$
whose Lie bracket is determined by
$$
([X,Y],Z)=\la Z(X),Y\ra, \qquad Z\in\tilde{\ngo}_2.
$$
In other words, the maps $\tilde{J}_Z$'s for this Lie algebra are
the $Z$'s themselves.  Recall that
$\dim{\tilde{\ngo}_2}=\frac{n(n-1)}{2}-k$, and so the dual
$\tilde{\ngo}$ of $\ngo$ is of type $(n,\frac{n(n-1)}{2}-k)$.  It is
proved in \cite{Sch} that $\ngo_1$ is isomorphic to $\ngo_2$ if and
only if $\tilde{\ngo}_1$ is isomorphic to $\tilde{\ngo}_2$, so that
any classification of type $(n,k)$ simultaneously classifies type
$(n,\frac{n(n-1)}{2}-k)$.

\begin{example}\label{dual42} {\rm
Let $\hg$ be the Lie algebra of type $(4,4)$ which is dual to
$\hg_3\oplus\hg_3$ (of type $(4,2)$).  The Lie bracket of
$\hg_3\oplus\hg_3$ is
$$
[X_1,X_2]=Z_1, \qquad [X_3,X_4]=Z_2,
$$
and hence
$$
J_{Z_1}=\left[\begin{smallmatrix} 0&-1&&\\ 1&0&&\\ &&0&0\\ &&0&0
\end{smallmatrix}\right],\qquad J_{Z_2}=\left[\begin{smallmatrix}
0&0&&\\ 0&0&&\\ &&0&-1\\ &&1&0 \end{smallmatrix}\right].
$$
The orthogonal complement $\tilde{\ngo}_2$ of $\{ J_Z:Z\in\ngo_2\}$
is then linearly generated by
$$
\left[\begin{smallmatrix} &&-1&0\\ &&0&0\\ 1&0&&\\ 0&0&&
\end{smallmatrix}\right], \quad \left[\begin{smallmatrix} &&0&-1\\
&&0&0\\ 0&0&&\\ 1&0&& \end{smallmatrix}\right], \quad
\left[\begin{smallmatrix} &&0&0\\ &&-1&0\\ 0&1&&\\ 0&0&&
\end{smallmatrix}\right], \quad \left[\begin{smallmatrix} &&0&0\\
&&0&-1\\ 0&0&&\\ 0&1&& \end{smallmatrix}\right],
$$
which determines the Lie bracket for $\hg$ given by}
\begin{equation}\label{corh}
[X_1,X_3]=Z_1, \quad [X_1,X_4]=Z_2, \quad [X_2,X_3]=Z_3, \quad
[X_2,X_4]=Z_4.
\end{equation}
\end{example}

Scheuneman duality also allows us to find all the rational forms of
$\hg$; namely, the duals of the rational form of $\hg_3\oplus\hg_3$,
already computed in Proposition \ref{rath3h3}.

\begin{proposition}\label{rath}
For any $k\in\ZZ$ let $\hg_k^{\QQ}$ be the rational Lie algebra of
type $(4,4)$ defined by
$$
\begin{array}{lcl}
[X_1,X_2]=Z_1, && [X_2,X_3]=-Z_3, \\

[X_1,X_3]=Z_2, && [X_2,X_4]=-Z_2, \\

[X_1,X_4]=kZ_3, && [X_3,X_4]=Z_4.
\end{array}
$$
Then the set of isomorphism classes of rational forms of the Lie
algebra $\hg$ defined in {\rm (\ref{corh})} is parametrized by
$$
\{ \hg_k^{\QQ}:k\,\mbox{is a square free natural number}\}.
$$
\end{proposition}

\begin{proof}
For the rational form $\ngoq_k$ of $\hg_3\oplus\hg_3$ (see
(\ref{crath3h3})) we have that
$$
J_{Z_1}=\left[\begin{smallmatrix} &&-1&0\\ &&0&-1\\ 1&0&&\\ 0&1&&
\end{smallmatrix}\right],\qquad J_{Z_2}=\left[\begin{smallmatrix}
&&0&-1\\ &&-k&0\\ 0&k&&\\ 1&0&& \end{smallmatrix}\right].
$$
A basis of the orthogonal complement of $\la
J_{Z_1},J_{Z_2}\ra_{\QQ}$ is then given by
$$
\left[\begin{smallmatrix} 0&-1&&\\ 1&0&&\\ &&0&0\\ &&0&0
\end{smallmatrix}\right], \quad \left[\begin{smallmatrix}&&-1&0\\
&&0&1\\ 1&0&&\\ 0&-1&& \end{smallmatrix}\right], \quad
\left[\begin{smallmatrix} &&0&-k\\ &&1&0\\ 0&-1&&\\ k&0&&
\end{smallmatrix}\right], \quad \left[\begin{smallmatrix} 0&0&&\\
0&0&&\\ &&0&-1\\ &&1&0
\end{smallmatrix}\right],
$$
which determines the Lie bracket for $\hg_k^{\QQ}$.  To conclude the
proof, one can easily check that the rational subspace generated by
$$
\begin{array}{l}
\left\{ \sqrt{k}(X_1-X_3), X_1+X_3, X_2+X_4, \sqrt{k}(X_2-X_4), \right. \\ \\

\left. 2\sqrt{k}Z_1, \sqrt{k}(Z_2+Z_3), Z_3-Z_2,
-2\sqrt{k}Z_4\right\},
\end{array}
$$
is closed under the Lie bracket of $\hg$ and isomorphic to
$\hg_k^{\QQ}$.
\end{proof}

An alternative proof of the non-isomorphism between the
$\hg^{\QQ}_k$'s without using Scheuneman duality may be given as
follows: from the form of $J_{Z_1},...,J_{Z_4}$ for $\hg^{\QQ}_k$ in
the above proof it follows that
$$
J_{xZ_1+yZ_2+zZ_3+wZ_4}=\left[\begin{smallmatrix} 0&-x&-y&-kz\\
x&0&z&y\\ y&-z&0&-w\\ kz&-y&w&0
\end{smallmatrix}\right],
$$
and so the Pfaffian form of $\hg^{\QQ}_k$ is given by
$f_k(x,y,z,w)=xw+y^2-kz^2$. Now, if $\hg^{\QQ}_k$ is isomorphic to
$\hg^{\QQ}_{k'}$ then $f_k\simeq_{\QQ} f_{k'}$ (see Proposition
\ref{isoforms}), which implies that $k=q^2k'$ for some $q\in\QQ^*$
by applying Proposition \ref{hessian} (recall that $Hf_k=4k$). Thus
$k=k'$ since they are square free.

We now study rational forms of $\lgo_4\oplus\lgo_4$, where $\lgo_4$
is the $4$-dimensional real Lie algebra with Lie bracket
$$
[Y_1,Y_2]=Y_3, \qquad [Y_1,Y_3]=Y_4.
$$
Since $\lgo_4\oplus\lgo_4$ is $3$-step nilpotent, Pfaffian forms and
duality can not be used as tools to distinguish or classify rational
forms, which makes of this case the hardest one.  For each
$k\in\ZZ$, consider the $8$-dimensional rational nilpotent Lie
algebra $\lgo_k^{\QQ}$ with basis $\{
X_1,X_2,X_3,X_4,Z_1,Z_2,Z_3,Z_4\}$ and Lie bracket defined by
\begin{equation}\label{deflk}
\begin{array}{lcl}
[X_1,X_3]=Z_1, && [X_2,X_3]=Z_2, \\

[X_1,X_4]=Z_2, && [X_2,X_4]=kZ_1, \\

[X_1,Z_1]=Z_3, && [X_2,Z_2]=kZ_3, \\

[X_1,Z_2]=Z_4, && [X_2,Z_1]=Z_4.
\end{array}
\end{equation}

{\small
\begin{table}
$$
\begin{array}{ccccc}
\hline\hline
&&& \\
Real\, Lie\, algebra & Type && Rational \; forms & Reference \\ \\
\hline
&&&\\
\hg_3\oplus\hg_3&(4,2) && \ngoq_k,\;k\ge 1& {\rm Prop.} \;\ref{rath3h3}  \\ \\

\fg_3           &(3,3) && \fg_3^\QQ       &   - - \\ \\

\ggo            &(6,2) && \ggo^\QQ & {\rm Prop.} \;\ref{rat62}  \\ \\

\hg_3\oplus\hg_5 & (6,2) && (\hg_3\oplus\hg_5)^{\QQ} & {\rm Prop.}\; \ref{rath3h5} \\ \\

\hg             &(4,4) && \hg^\QQ_k, \; k\ge 1& {\rm Prop.} \;\ref{rath}   \\ \\

\lgo_4\oplus\lgo_4 & (4,2,2) && \lgo_k^\QQ, \; k\ge 1& {\rm Prop.}\; \ref{ratl4l4} \\ \\

\hline \hline \\
\end{array}
$$
\caption{Set of rational forms up to isomorphism for some real
nilpotent Lie algebras. In all cases $k$ runs over all square-free
natural numbers.}\label{ratforms}
\end{table}}

\begin{theorem}\label{ratl4l4}
Let $\{ X_1,X_2,X_3,X_4,Z_1,Z_2,Z_3,Z_4\}$ be a basis of the Lie
algebra $\lgo_4\oplus\lgo_4$ of type $(4,2,2)$ with structure
coefficients
$$
\begin{array}{lcl}
[X_1,X_3]=Z_1, && [X_2,X_4]=Z_2, \\

[X_1,Z_1]=Z_3, && [X_2,Z_2]=Z_4.
\end{array}
$$
For each $k\in\NN$ the rational subspace generated by the set
$$
\begin{array}{l}
\left\{ X_1+X_2, \sqrt{k}(X_1-X_2), X_3+X_4, \sqrt{k}(X_3-X_4), \right. \\

\left. Z_1+Z_2, \sqrt{k}(Z_1-Z_2), Z_3+Z_4,
\sqrt{k}(Z_3-Z_4)\right\}
\end{array}
$$
is a rational form of $\lgo_4\oplus\lgo_4$ isomorphic to the Lie
algebra $\lgo_k^{\QQ}$ defined in {\rm (\ref{deflk})}.  Moreover,
the set
$$
\{\lgo_k^{\QQ}:k\;\mbox{is a square-free natural number}\}
$$
parametrizes all the rational forms of $\lgo_4\oplus\lgo_4$ up to
isomorphism.
\end{theorem}

\begin{proof}
It is easy to see that the Lie brackets of the basis of the rational
subspace coincides with the one of $\lgo_k^{\QQ}$ by renaming the
basis as $\{ X_1,...,Z_4\}$ with the same order.  In particular,
such a subspace is a rational form of $\lgo_4\oplus\lgo_4$.  If
$k'=q^2k$ then one can easily check that
$A:\lgo_{k'}^{\QQ}\mapsto\lgo_k^{\QQ}$ given by the diagonal matrix
with entries $(1,q,1,q,1,q,1,q)$ is an isomorphism of Lie algebras.

Conversely, assume that $A:\lgo_{k}^{\QQ}\mapsto\lgo_{k'}^{\QQ}$ is
an isomorphism. We will show that $k'=q^2k$ for some $q\in\QQ^*$.
Let $\{ J_Z'\}$, $\{ J_Z\}$ be the maps defined at the beginning of
this section corresponding to $\lgo_{k'}^{\QQ}$ and $\lgo_k^{\QQ}$,
respectively.  If $Z=xZ_1+yZ_2+zZ_3+wZ_4$ we have that
$$
J_Z=\left[\begin{array}{cccccccc}
0&0&-x&-y&-z&-w&0&0\\
0&0&-y&-kx&-w&-kz&0&0\\
x&y&0&&\cdots&&&0\\
y&kx&&&&&&\\
z&w&\vdots&&&&\vdots&\vdots\\
w&kz&&&&&&\\
0&0&&&&&&\\
0&0&0&&\cdots&&&0
\end{array}\right],
$$
and $J_Z'$ is obtained just by replacing $k$ with $k'$.  It follows
from (\ref{iso}) that $A^tJ_Z'A=J_{A^tZ}$ for all $Z\in\la
Z_3,Z_4\ra_{\QQ}$, and since this subspace is $A$-invariant we get
that the subspace
$$
\displaystyle{\bigcap_{Z\in\la Z_3,Z_4\ra_{\QQ}}} \Ker{J_Z}=
\displaystyle{\bigcap_{Z\in\la Z_3,Z_4\ra_{\QQ}}} \Ker{J'_Z} =\la
X_3,X_4,Z_3,Z_4\ra_{\QQ}
$$
is also $A$-invariant.  Thus $A$ has the form
\begin{equation}\label{Aform}
A=\left[\begin{array}{cccc}
A_1&0&0&0\\
\star&A_2&0&0\\
\star&0&A_3&0\\
\star&\star&\star&A_4.
\end{array}\right]
\end{equation}
(recall that $C^1(\lgo_k^{\QQ})=C^1(\lgo_{k'}^{\QQ})=\la
Z_1,Z_2,Z_3,Z_4\ra_{\QQ}$ and
$C^2(\lgo_k^{\QQ})=C^2(\lgo_{k'}^{\QQ})=\la Z_3,Z_4\ra_{\QQ}$ are
always $A$-invariant), and now it is easy to prove that
$$
A_3^t\left[\begin{array}{cc}z&w\\ w&k'z\end{array}\right]A_1=  \left[\begin{array}{cc}az+bw&cz+dw\\
cz+dw&k'(az+bw)\end{array}\right], \; {\rm where}\;
A^t_4=\left[\begin{array}{cc}a&bw\\ c&d\end{array}\right].
$$
We compute the determinant of both sides getting
$$
qf'(z,w)=f(A_4^t(z,w)), \qquad \forall\; (z,w)\in\QQ^2,
$$
where $q=\det{A_3A_1}\in\QQ^*$ and $f(z,w)=kz^2-w^2$,
$f'(z,w)=k'z^2-w^2$.  By Proposition \ref{hessian} we have that
$$
4k'=q^{-2}(\det{A_4})^24k,
$$
and so $k=k'$ as long as they are square free numbers, as we wanted
to show.

To conclude the proof, it remains to show that these are all the
rational forms up to isomorphism.  Let $\ngoq$ be a rational form of
$\lgo_4\oplus\lgo_4$.  Since $\ngoq/[\ngoq,[\ngoq,\ngoq]]$ is of
type $(4,2)$, we can use the classification of rational Lie algebras
of this type given in (\ref{crath3h3}) to get linearly independent
vectors $X_1,...,Z_2$ such that
\begin{equation}\label{crat}
[X_1,X_3]=Z_1, \quad [X_1,X_4]=Z_2, \quad [X_2,X_3]=Z_2, \quad
[X_2,X_4]=kZ_1,
\end{equation}
where $k$ is a square free integer number.  Jacobi condition is
equivalent to
\begin{equation}\label{jacobi}
\begin{array}{lcl}
[X_1,Z_2]=[X_2,Z_1], && [X_3,Z_2]=[X_4,Z_1], \\ \\
k[X_1,Z_1]=[X_2,Z_2], && k[X_3,Z_1]=[X_4,Z_2].
\end{array}
\end{equation}
We will consider the following two cases separately:
\begin{itemize}
\item[(I)]  $Z_3:=[X_1,Z_1]$ and $Z_4:=[X_1,Z_2]$ are linearly independent,

\item[(II)] $[X_1,Z_1], [X_1,Z_2] \in\QQ Z_3$ for some nonzero $Z_3\in\ngoq$.
\end{itemize}
In both cases we will make use of the following isomorphism
invariant for real $3$-step nilpotent Lie algebras:
$$
U(\ngo):=\left\{
X\in\ngo/[\ngo,[\ngo,\ngo]]:\dim{\Im(\ad{X})}=1\right\}\cup\{ 0\}.
$$
Clearly, if $A:\ngo\mapsto\ngo'$ is an isomorphism then
$AU(\ngo)=U(\ngo')$.  Under the presentation of $\lgo_4\oplus\lgo_4$
given in the statement of the theorem, it is easy to see that
\begin{equation}\label{ul4l4}
U(\lgo_4\oplus\lgo_4)=\la X_3,Z_1\ra_{\RR}\cup\la X_4,Z_3\ra_{\RR}.
\end{equation}
In case (I), it follows from (\ref{jacobi}) that we also have
$$
[X_2,Z_1]=Z_4, \qquad [X_2,Z_2]=kZ_3.
$$
Therefore, in order to get that $\ngoq$ is isomorphic to
$\lgo_k^{\QQ}$ (see (\ref{deflk})), it is enough to show that the
vectors in $\la Z_3,Z_4\ra_{\RR}$ given by
$$
Z:=k[X_3,Z_1]=[X_4,Z_2], \qquad Z':=[X_3,Z_2]=[X_4,Z_1]
$$
are both zero (see (\ref{jacobi})).  Let us compute the cone
$U(\ngo)$ for $\ngo=\ngoq\otimes\RR$.  Recall that $U(\ngo)$ has to
be the union of two disjoint planes as
$\ngo\simeq\lgo_4\oplus\lgo_4$ (see (\ref{ul4l4})).  If
$X=aX_1+bX_2+cX_3+dX_4+eZ_1+fZ_2$ then
$$
\begin{array}{l}
[X_1,X]=cZ_1+dZ_2+eZ_3+fZ_4, \\

[X_2,X]=dkZ_1+cZ_2+fkZ_3+eZ_4, \\

[X_3,X]=-aZ_1-bZ_2+\frac{e}{k}Z+fZ', \\

[X_4,X]=-bkZ_1-aZ_2+fZ+eZ', \\

[Z_1,X]=-aZ_3-bZ_4-\frac{c}{k}Z-d Z', \\

[Z_2,X]=-bkZ_3-aZ_4-d Z-c Z'.
\end{array}
$$

Assume that $\Im(\ad{X})=\RR X_0$, $X_0\ne 0$.  If $k\leq 0$ then it
follows easily from $[X_1,X]=\l[X_2,X]$ and $[X_3,X]=\mu[X_4,X]$ for
some $\l,\mu\in\RR$ that $a=b=c=d=e=f=0$, which implies that
$U(\ngo)=\{ 0\}$, a contradiction.

\begin{remark}\label{frl4l4}
{\rm Since $k$ has to be positive one can also get by an easy
adaptation of this proof that the only real form of
$(\lgo_4\oplus\lgo_4)_{\CC}$ is $\lgo_4\oplus\lgo_4$.}
\end{remark}

We then have that $k>0$ and $a=\pm\sqrt{k}b$, $c=\pm\sqrt{k}d$,
$e=\pm\sqrt{k}f$, where $c$ and $e$ have the same sign.  This
implies that
$$X=b(\pm\sqrt{k}X_1+X_2)+d(\pm\sqrt{k}X_3+X_4)+f(\pm\sqrt{k}Z_1+Z_2)$$
and
$$
\begin{array}{l}
[X_1,X]=d(\pm\sqrt{k}Z_1+Z_2)+f(\pm\sqrt{k}Z_3+Z_4), \\

[X_2,X]=\sqrt{k}[X_1,X], \\

[X_3,X]=-b(\pm\sqrt{k}Z_1+Z_2)+f(\pm\frac{1}{\sqrt{k}}Z+Z'), \\

[X_4,X]=\sqrt{k}[X_3,X], \\

[Z_1,X]=-b(\pm\sqrt{k}Z_3+Z_4)-d(\pm\frac{1}{\sqrt{k}}Z+Z'), \\

[Z_2,X]=\sqrt{k}[Z_1,X].
\end{array}
$$

If $b\ne 0$ then $d\ne 0$ and $a$ has the same sign as $c$ and $e$,
and since $X_0$ has a nonzero component in $\la Z_1,Z_2\ra_{\RR}$ we
get $[Z_1,X]=0$, that is,
$-\frac{b}{d}(\pm\sqrt{k}Z_3+Z_4)=\pm\frac{1}{\sqrt{k}}Z+Z'$.  In
any case we obtain a subset of $U(\ngo)$ of the form
$$
\{ b(\pm\sqrt{k}X_1+X_2)+
d(\pm\sqrt{k}X_3+X_4)+f(\pm\sqrt{k}Z_1+Z_2):b,d\ne 0\}
$$
with the same sign in all the terms, which is a contradiction since
$U(\ngo)$ is the union of two planes.  Thus $b=0$ and so
$$
U(\ngo)=\la\sqrt{k}X_3+X_4,
\sqrt{k}Z_1+Z_2\ra_{\RR}\cup\la-\sqrt{k}X_3+X_4,
-\sqrt{k}Z_1+Z_2\ra_{\RR}.
$$
This clearly implies that
$\frac{1}{\sqrt{k}}Z+Z'=-\frac{1}{\sqrt{k}}Z+Z'=0$, that is
$Z=Z'=0$, as was to be shown.

Concerning case (II), we can assume that
$$
[X_1,Z_2]=rZ_3, \quad k[X_1,Z_1]=sZ_3, \quad [X_3,Z_2]=tZ_4, \quad
k[X_3,Z_1]=uZ_4,
$$
where $Z_3,Z_4$ are linearly independent and $(s,r),(u,t)\ne (0,0)$.
By using (\ref{jacobi}), for $X=aX_1+bX_2+cX_3+dX_4+eZ_1+fZ_2$ we
have that
$$
\begin{array}{l}
[X_1,X]=cZ_1+dZ_2+(\frac{e}{k}s+fr)Z_3, \\

[X_2,X]=dkZ_1+cZ_2+(fs+er)Z_3, \\

[X_3,X]=-aZ_1-bZ_2+(\frac{e}{k}u+ft)Z_4, \\

[X_4,X]=-bkZ_1-aZ_2+(fu+et)Z_4, \\

[Z_1,X]=-(\frac{a}{k}s+br)Z_3-(\frac{c}{k}u+dt)Z_4, \\

[Z_2,X]=-(\frac{b}{k}s+ar)Z_3-(\frac{d}{k}u+ct)Z_4.
\end{array}
$$

If $a=0$ then $b=c=d=0$.  We also obtain that $e^2=kf^2$, since
either
$$
\left[\begin{array}{cc} \frac{e}{k} & f \\ f & e \end{array}\right] \left[\begin{array}{c} s \\
r \end{array}\right]=0 \qquad \mbox{or} \qquad
\left[\begin{array}{cc} \frac{e}{k} & f \\ f & e \end{array}\right] \left[\begin{array}{c} u \\
t \end{array}\right]=0.
$$
We do not get any plane in $U(\ngo)$ in this way and therefore there
must be an $X\in U(\ngo)$ with $a\ne 0$, which implies that
$b,c,d\ne 0$ and $a^2=kb^2$, $c^2=kd^2$. Thus $[Z_1,X]= [Z_2,X]=0$
and so $\Im(\ad{X})\subset \la Z_1,Z_2\ra_\RR.$ This implies that
$e^2=kf^2$ and then the $3$-dimensional subspace
$$
\la\sqrt{k}X_1+X_2, \sqrt{k}X_3+X_4, \sqrt{k}Z_1+Z_2\ra_{\RR}\subset
U(\ngo),
$$
which is a contradiction, proving that case (II) is not possible and
concluding the proof.

\end{proof}

\section{Anosov diffeomorphisms and Lie algebras}\label{andiff}

Anosov diffeomorphisms play an important and beautiful role in
dynamics as the notion represents the most perfect kind of global
hyperbolic behavior, giving examples of structurally stable
dynamical systems.  A diffeomorphism $f$ of a compact differentiable
manifold $M$ is called {\it Anosov} if the tangent bundle $\tang M$
admits a continuous invariant splitting $\tang M=E^+\oplus E^-$ such
that $\dif f$ expands $E^+$ and contracts $E^-$ exponentially, that
is, there exist constants $0<c$ and $0<\lambda<1$ such that
$$
||\dif f^n(X)||\leq c\lambda^n||X||, \quad \forall X\in E^-, \qquad
||\dif f^n(Y)||\geq c\lambda^{-n}||Y||, \quad \forall Y\in E^+,
$$
for all $n\in\NN$.  The condition is independent of the Riemannian
metric.  Some of the other very nice properties of these special
dynamical systems, all proved mainly by D. Anosov, are: the
distributions $E^+$ and $E^-$ are completely integrable with
$C^{\infty}$ leaves and determine two (unique) $f$-invariant
foliations (unstable and stable, respectively) with remarkable
dynamical properties; the set of periodic points (i.e. $f^m(p)=p$
for some $m\in\NN$) is dense in the set of those points of $M$ such
that for any neighborhood $U$ of $p$ there exist $k\ne m\in\NN$ with
$f^k(U)\cap f^m(U)\ne\emptyset$; the set of all Anosov
diffeomorphisms form an open subset of $\Diff(M)$.

\begin{example}\label{nilmanifolds} {\rm Let $N$ be a real simply connected nilpotent
Lie group with Lie algebra $\ngo$. Let $\vp$ be a hyperbolic
automorphism of $N$, that is, all the eigenvalues of its derivative
$A=(\dif\vp)_e:\ngo\mapsto\ngo$ have absolute value different from
$1$. If $\vp(\Gamma)=\Gamma$ for some lattice $\Gamma$ of $N$ (i.e.
a uniform discrete subgroup) then $\vp$  defines an Anosov
diffeomorphism on the nilmanifold $M=N/\Gamma$, which is called an
{\it Anosov automorphism}.  The subspaces $E^+$ and $E^-$ are
obtained by left translation of the eigenspaces of eigenvalues of
$A$ of absolute value greater than $1$ and less than $1$,
respectively, and so the splitting is differentiable.  If more in
general, $\Gamma$ is a uniform discrete subgroup of $K\ltimes N$,
where $K$ is any compact subgroup of $\Aut(N)$, for which
$\vp(\Gamma)=\Gamma$ (recall that $\vp$ acts on $\Aut(N)$ by
conjugation), then $\vp$ also determines an Anosov diffeomorphism on
$M=N/\Gamma$ which is also called Anosov automorphism.  In this case
$M$ is called an infranilmanifold and is finitely covered by the
nilmanifold $N/(N\cap\Gamma)$. }
\end{example}

In \cite{Sml}, S. Smale raised the problem of classifying all
compact manifolds (up to homeomorphism) which admit an Anosov
diffeomorphism.  At this moment, the only known examples are of
algebraic nature, namely Anosov automorphisms of nilmanifolds and
infranilmanifolds described in the example above.  It is conjectured
that any Anosov diffeomorphism is topologically conjugate to an
Anosov automorphism of an infranilmanifold (see \cite{Mrg}).

All this certainly highlights the problem of classifying all
nilmanifolds which admit Anosov automorphisms, which are easily seen
in correspondence with a very special class of nilpotent Lie
algebras over $\QQ$.  Nevertheless, not too much is known on the
question since it is not so easy for an automorphism of a (real)
nilpotent Lie algebra being hyperbolic and unimodular at the same
time.

\begin{definition}\label{anolie}
{\rm A rational Lie algebra $\ngoq$ (i.e. with structure constants
in $\QQ$) of dimension $n$ is said to be {\it Anosov} if it admits a
{\it hyperbolic} automorphism $A$ (i.e. all their eigenvalues have
absolute value different from $1$) which is {\it unimodular}, that
is, $[A]_{\beta}\in\Gl_n(\ZZ)$ for some basis $\beta$ of $\ngoq$,
where $[A]_{\beta}$ denotes the matrix of $A$ with respect to
$\beta$.  A hyperbolic and unimodular automorphism is called an {\it
Anosov automorphism}.  We also say that a real Lie algebra is {\it
Anosov} when it admits a rational form which is Anosov.  An
automorphism of a real Lie algebra $\ngo$ is called {\it Anosov} if
it is hyperbolic and $[A]_{\beta}\in\Gl_n(\ZZ)$ for some $\ZZ$-basis
$\beta$ of $\ngo$ (i.e. with integer structure constants).}
\end{definition}

The unimodularity condition on $A$ in the above definition is
equivalent to the fact that the characteristic polynomial of $A$ has
integer coefficients and constant term equal to $\pm 1$ (see
\cite{Dkm}).  It is well known that any Anosov Lie algebra is
necessarily nilpotent, and it is easy to see that the classification
of nilmanifolds which admit an Anosov automorphism is essentially
equivalent to that of Anosov Lie algebras (see
\cite{anosov,Dn,Ito,Dkm}). If $\ngo$ is a rational Lie algebra, we
call the real Lie algebra $\ngo \otimes \RR$ the {\it real
completion} of $\ngo$.

We now give some necessary conditions a real Lie algebra has to
satisfy in order to be Anosov (see \cite{Mlf}).

\begin{proposition}\label{auto}
Let $\ngo$ be a real nilpotent Lie algebra which is Anosov.  Then
there exist a decomposition $\ngo=\ngo_1\oplus...\oplus\ngo_r$
satisfying $C^i(\ngo)=\ngo_{i+1}\oplus...\oplus\ngo_r$, $i=0,...,r$,
and a hyperbolic $A\in\Aut(\ngo)$ such that
\begin{itemize}
\item[(i)] $A\ngo_i=\ngo_i$ for all $i=1,...,r$.

\item[(ii)] $A$ is semisimple (in particular $A$ is diagonalizable over $\CC$).

\item[(iii)] For each $i$, there exists a
basis $\beta_i$ of $\ngo_i$ such that
$[A_i]_{\beta_i}\in\Sl_{n_i}(\ZZ),$ where $n_i=\dim{\ngo_i}$ and
$A_i=A|_{\ngo_i}$.
\end{itemize}
\end{proposition}

\begin{proof}
Let $\beta$ be a $\ZZ$-basis of $\ngo$ for which there is a
hyperbolic $A\in\Aut(\ngo)$ satisfying $[A]_{\beta}\in\Gl_n(\ZZ)$.
By using that $\Aut(\ngo)$ is a linear algebraic group, it is proved
in \cite[Section 2]{AslSch} that we can assume that $A$ is
semisimple.  Thus the existence of the decomposition satisfying (i)
follows from the fact that the subspaces $C^i(\ngo)$ are
$A$-invariant.

If $\beta=\{ X_1,...,X_n\}$ then the discrete (additive) subgroup
$$
\ngoz=\left\{\sum_{i=1}^na_iX_i:a_i\in\ZZ\right\}
$$
of $\ngo$ is closed under the Lie bracket of $\ngo$ and
$A$-invariant, and $C^i(\ngoz)$ is a discrete subgroup of
$C^i(\ngo)$ of maximal rank.  Since $AC^i(\ngoz)=C^i(\ngoz)$ for any
$i$ we have that $A$ induces an invertible map
$$
C^{i-1}(\ngoz)/C^i(\ngoz)\mapsto C^{i-1}(\ngoz)/C^i(\ngoz),
$$
and it follows from $C^i(\ngoz)\otimes\RR=C^i(\ngo)$ that
$C^{i-1}(\ngoz)/C^i(\ngoz)\simeq\ZZ^{n_i}$ is a discrete subgroup of
$C^{i-1}(\ngo)/C^i(\ngo)\simeq\ngo_i$ which is leaved invariant by
$A$, proving the existence of the basis $\beta_i$ of $\ngo_i$ in
(iii).  Recall that by considering $A^2$ rather than $A$ if
necessary, we can assume that $\det{A_i}=1$ for all $i$.
\end{proof}

\begin{proposition}\label{coroauto}
Let $\ngo$ be a real $r$-step nilpotent Lie algebra of type
$(n_1,...,n_r)$ (see {\rm Definition \ref{type}}). If $\ngo$ is
Anosov then at least one of the following is true:
\begin{itemize}
\item[(i)] $n_1\geq 4$ and $n_i\geq 2$ for all $i=2,...,r$.

\item[(ii)] $n_1=n_2=3$ and $n_i\geq 2$ for all $i=3,...,r$.
\end{itemize}
In particular, $\dim{\ngo}\geq 2r+2$.
\end{proposition}

\begin{proof}
We know from Proposition \ref{auto} that $A_i\in\Sl_{n_i}(\ZZ)$ is
hyperbolic, which implies that $n_i\geq 2$ for any $i$.  Assuming
(i) does not hold means then that $n_1=3$.  If $n_2=2$ and
$\{\lambda_1,\lambda_2,\lambda_3\}$ are the eigenvalues of $A_1$
then the eigenvalues of $A_2$ are of the form $\lambda_i\lambda_j$,
say $\{\lambda_1\lambda_2,\lambda_1\lambda_3\}$, and hence
$\lambda_1=\lambda_1^2\lambda_2\lambda_3=1$, which contradicts the
fact that $A_1$ is hyperbolic.  This implies that $n_2=3$.
\end{proof}

In \cite[Question (ii)]{anosov} there are examples of real Anosov
Lie algebras of type $(4,2,...,2)$ for any $r\geq 2$.  We shall
prove in Section \ref{cla}, Case $(3,3,2)$, that in part (ii) of the
above proposition one actually needs $n_3\geq 3$. Also, we do not
know of any example of type of the form $(3,3,...)$.

Part (i) of the following proposition is essentially \cite[Theorem
3]{AslSch}

\begin{proposition}\label{region}
Let $\ngo=\ngo_1\oplus\ngo_2$ be a real $2$-step nilpotent Lie
algebra with $\dim{\ngo_2}=k$.  Assume that $\ngo$ is Anosov and let
$\ngoq$ denote the rational form which is Anosov.
\begin{itemize}
\item[(i)] If $f$ is the Pfaffian form of $\ngo$ then for any $c>0$ the region
$$
R_c=\{(x_1,...,x_k)\in\RR^k:|f(x_1,...,x_k)|\leq c\}
$$
is unbounded.

\item[(ii)] For the Pfaffian form $f$ of $\ngoq$ and for any $p\in\ZZ$ the set
$$
S_p=\{(x_1,...,x_k)\in\ZZ^k:f(x_1,...,x_k)=p\}
$$
is either empty or infinite.
\end{itemize}
\end{proposition}

\begin{proof}
(i) Consider $A\in\Aut(\ngo)$ satisfying all the conditions in
Proposition \ref{auto}.  It follows from the proof of Proposition
\ref{isoforms} and $\det{A_i}=1$ for any $i=1,...,r$ that
$$
f(x_1,...,x_k)=f(A^t(x_1,...,x_k)) \qquad \forall\;
(x_1,...,x_k)\in\RR^k=\ngo_2,
$$
and so $AR_c\subset R_c$.  Assume that $R_{c_0}$ is bounded for some
$c_0>0$, by using that $f$ is a homogeneous polynomial we get that
$R_c$ is bounded for any $c>0$; indeed, $R_c=c^{-\frac{1}{m}}R_1$ if
$m$ is the degree of $f$.  Now, for a sufficiently big $c_1>0$ we
may assume that $R_{c_1}$ contains the basis $\beta_2$ of $\ngo_2$,
but only finitely many integral linear combinations of elements in
this basis can belong to the bounded region $R_{c_1}$.  This implies
that $A^t|_{\ngo_2}$ leave a finite set of points invariant, and
since such a set contains a basis of $\ngo_2$ we obtain that
$(A^t)^l=I$ for some $l\in\NN$. The eigenvalues of $A$ have then to
be roots of the identity, contradicting its hyperbolicity.

\no (ii) Analogously to the proof of part (i), we get that
$A^tS_p\subset S_p$.  If $S_p\ne\emptyset$ and finite then for the
real subspace $W\subset \ngo_2$ generated by $S_p$ we have that
$A^tW\subset W$ and $(A^t|_W)^l=I$ for some $l\in\NN$, which is
again a contradiction by the hyperbolicity of $A$.
\end{proof}

We now give an example of how the above proposition can be applied.
Rational Lie algebras of type $(4,2)$ are parametrized by the set of
square free numbers $k \in \ZZ$ and their Pfaffian forms are
$f_k(x,y) = x^2-k y^2$ (see the paragraph before Proposition
\ref{rath3h3}).  Thus the set of solutions
$$ \left\{ (x,y) \in \ZZ^2 \;:\; f_k(x,y)=1 \right\} $$
is infinite if and only if $k>1$ or $k=0$ (Pell equation).  By
Proposition \ref{region}, (ii), the Lie algebra $\ngoq_k$ can never
be Anosov for $k<0$ or $k=1$.  Recall that we could also discard
$\ngoq_k,$ $k<0$ as a real Anosov Lie algebra by applying
Proposition \ref{region}, (i).

It is not true in general that if a direct sum of real Lie algebras
is Anosov then each of the direct summands is so, as the example
$\hg_3\oplus\hg_3$ shows. However, we shall see next that this
actually happens when one of the direct summands is (maximal)
abelian.

Let $\ngo$ be a Lie algebra over $K$. An {\it abelian factor} of
$\ngo$ is an abelian ideal $\ag$ for which there exists an ideal
$\tilde{\ngo}$ of $\ngo$ such that $\ngo=\tilde{\ngo}\oplus\ag$
(i.e. $[\tilde{\ngo},\ag]=0).$  Let $m(\ngo)$ denote the maximum
dimension over all abelian factors of $\ngo$.  If $\zg$ is the
center of $\ngo$ then the maximal abelian factors are precisely the
linear direct complements of $\zg \cap [\ngo,\ngo]$ in $\zg,$ that
is, those subspaces $\ag \subset \zg$ such that $\zg = \zg \cap
[\ngo,\ngo] \oplus \ag.$  Therefore
$$ m(\ngo)= \dim \zg - \dim \zg \cap [\ngo,\ngo].$$

\begin{theorem}\label{abfactor}
Let $\ngo$ be a rational Lie algebra with $m(\ngo)=r$ and let
$\ngo=\tilde{\ngo}\oplus\QQ^r$ be any decomposition in ideals, that
is, $\QQ^r$ is a maximal abelian factor of $\ngo.$  Then $\ngo$ is
Anosov if and only if $\tilde{\ngo}$ is Anosov and $r \ge 2.$
\end{theorem}

\begin{proof}
If $\tilde{\ngo}$ is Anosov and $r \ge 2$ then we consider the
automorphism $A$ of $\ngo$ defined on $\tilde{\ngo}$ as an Anosov
automorphism of $\tilde{\ngo}$ and on $\QQ^r$ as any hyperbolic
matrix in $\Gl_r(\ZZ).$  Thus $A$ is an Anosov automorphism of
$\ngo.$

Conversely, let $A$ be an Anosov automorphism of $\ngo.$  As in the
proof of Proposition \ref{auto} we may assume that $A$ is semisimple
and consider the discrete (additive) subgroup
$$ \ngo^\ZZ = \left\{ \sum^n_{i=1} a_i X_i, \; a_i \in \ZZ \right\}$$
which is $A$-invariant.  Since the center $\zg$ of $\ngo$ and $\zg_1
= \zg \cap [\ngo,\ngo]$ are both leaved invariant by $A$, there
exist $A$-invariant subspaces $V$ and $\ag \subset \zg$ such that
$$
\ngo = V \oplus \zg = V \oplus \zg_1 \oplus \ag.
$$
Thus $\ag$ is a maximal abelian factor, $\dim\ag = r$ and $A$ has
the form
$$
\begin{array}{llll}
 A= \left[\begin{smallmatrix}A_1&&\\&A_2&\\&&A_3\end{smallmatrix}\right],& A_1=A|_{V}, &
A_2=A|_{\zg_1},& A_3=A|_{\ag}.
\end{array}
$$
The subgroup $\zg(\ngo_\ZZ) =\left\{ X \in \ngo_\ZZ : [X,Y]=0 \;\;
\forall \; Y \in \ngo_\ZZ \right\}$ is also $A$-invariant and it is
a lattice of $\zg$ (i.e. a discrete subgroup of maximal rank) since
for any $Z \in \zg$ there exist $k \in \ZZ$ such that $kZ \in
\zg(\ngo_\ZZ)$ and $Z=\frac{1}{k}(kZ),$ that is,
$\zg(\ngo_\ZZ)\otimes \QQ = \zg.$ Since $\ngo_\ZZ / \zg(\ngo_\ZZ)$
is $A$-invariant and $\left(\ngo_\ZZ / \zg(\ngo_\ZZ)\right) \otimes
\QQ \simeq V $ we get that $A_1$ is unimodular. Analogously, $A_2$
and $A_3$ are unimodular since $\zg_1(\ZZ)=\zg(\ngo_\ZZ)\cap
[\ngo_\ZZ,\ngo_\ZZ]$ and $\zg(\ngo_\ZZ)/\zg_1(\ZZ)$ are also
discrete subgroups of maximal rank of $\zg_1$ and $\zg/\zg_1 \simeq
\ag,$ respectively.

The hyperbolicity of $A$ guaranties the one of $A_1, A_2$ and $A_3$
and so $\tilde{\ngo} \simeq V \oplus \zg_1$ is Anosov and $\dim \ag
\ge 2$, as we wanted to show.
\end{proof}

To finish this section, we give a simple procedure to construct
explicit examples of Anosov Lie algebras.  This result is a
generalization of \cite[Theorem 3.1]{anosov} proposed by F.
Grunewald.

A Lie algebra $\ngo$ over $K$ is said to be {\it graded} (over
$\NN$) if there exist $K$-subspaces $\ngo_i$ of $\ngo$ such that
$$
\ngo=\ngo_1\oplus\ngo_2\oplus...\oplus\ngo_k \qquad \mbox{and}
\qquad [\ngo_i,\ngo_j]\subset\ngo_{i+j}.
$$
Equivalently, $\ngo$ is graded when there are nonzero $K$-subspaces
$\ngo_{d_1},...,\ngo_{d_r}$, $d_1<...<d_r$, such that
$\ngo=\ngo_{d_1}\oplus...\oplus\ngo_{d_r}$ and if $0\ne
[\ngo_{d_i},\ngo_{d_j}]$ then $d_i+d_j=d_k$ for some $k$ and
$[\ngo_{d_i},\ngo_{d_j}]\subset \ngo_{d_k}$. Recall that any graded
Lie algebra is necessarily nilpotent.

\begin{theorem}\label{gradedsum}
Let $\ngoq$ be a graded rational Lie algebra, and consider the
direct sum $\tilde{\ngo}^{\QQ}=\ngoq\oplus...\oplus\ngoq$ ($s$
times, $s\geq 2$).  Then the real Lie algebra
$\tilde{\ngo}=\tilde{\ngo}^{\QQ}\otimes\RR$ is Anosov.  In other
words, if $\ngo$ is a graded real Lie algebra admitting a rational
form, then $\tilde{\ngo}=\ngo\oplus...\oplus\ngo$ ($s$-times, $s\geq
2$)  is Anosov.
\end{theorem}

\begin{proof}
Let $\{ X_1,...,X_n\}$ be a $\ZZ$-basis of $\ngoq$ compatible with
the gradation $\ngoq=\ngoq_{d_1}\oplus...\oplus\ngoq_{d_r}$, that
is, a basis with integer structure constants and such that each
$X_i\in\ngoq_{d_j}$ for some $j$. We will denote this basis by $\{
X_{l1},...,X_{ln}\}$ when we need to make clear that it is a basis
of the $l$-th copy of $\ngoq$ in $\tilde{\ngo}^{\QQ}$, so the Lie
bracket of $\tilde{\ngo}^{\QQ}$ is given by $[X_{li},X_{l'j}]=0$ for
all $l\ne l'$, and for any $l=1,...,s$
\begin{equation}\label{corchete}
[X_{li},X_{lj}]=\sum_{k=1}^{n} m_{ij}^kX_{lk}, \qquad
m_{ij}^k\in\ZZ.
\end{equation}
Every nonzero $\lambda\in\RR$ defines an automorphism $A_{\lambda}$
of $\ngoq\otimes\RR$ by
$$
A_{\lambda}|_{\ngoq_{d_i}\otimes\RR}=\lambda^{d_i}I.
$$
Let $B$ be a matrix in $\Gl_s(\ZZ)$ with eigenvalues
$\lambda_1,...,\lambda_s$ and assume that all of them are real
numbers different from $\pm 1$ (we are using here that $s\geq 2$).
This determines an automorphism $A$ of $\tilde{\ngo}$ in the
following way: $A$ leaves the decomposition
$\tilde{\ngo}^{\QQ}=(\ngoq\otimes\RR)\oplus...\oplus(\ngoq\otimes\RR)$
invariant and on the $l$-th copy of $\ngoq\otimes\RR$ coincides with
$A_{\lambda_l}$.

Consider the new basis of $\tilde{\ngo}$ defined by
$$
\begin{array}{rl}
\beta=&\{ X_{11}+X_{21}+...+X_{s1},
\lambda_1X_{11}+\lambda_2X_{21}+...+\lambda_sX_{s1},..., \\ \\
&\lambda_1^{s-1}X_{11}+\lambda_2^{s-1}X_{21}+...
+\lambda_s^{s-1}X_{s1},...,...,X_{1n}+X_{2n}+...+X_{sn}, \\ \\
&\lambda_1X_{1n}+\lambda_2X_{2n}+...+\lambda_sX_{sn},...,
\lambda_1^{s-1}X_{1n}+\lambda_2^{s-1}X_{2n}+...+\lambda_s^{s-1}X_{sn}
\}.
\end{array}
$$
In order to prove that $\beta$ is also a $\ZZ$-basis we take two
generic elements of it, say
$X=\lambda_1^{t}X_{1i}+\lambda_2^{t}X_{2i}+...+\lambda_s^tX_{si}$
and
$Y=\lambda_1^{u}X_{1j}+\lambda_2^{u}X_{2j}+...+\lambda_s^uX_{sj}$
for some $0\leq t,u\leq s-1$ and $1\leq i,j\leq n$.  Since the
$\lambda_l$'s are all roots of the characteristic polynomial
$f(x)=a_0+a_1x+...+a_{s-1}x^{s-1}+x^s$ of $B$ (with $a_i\in\ZZ$ and
$a_0=\pm 1$), there exist $b_0,...,b_{s-1}\in\ZZ$ (independent from
$l$) such that
$\lambda_l^{t+u}=b_0+b_1\lambda_l+...+b_{s-1}\lambda_l^{s-1}$ for
any $l=1,...,s$.  Now, by using (\ref{corchete}) we obtain that
$$
\begin{array}{rl}
[X,Y]&=\lambda_1^{t+u}[X_{1i},X_{1j}]+...+\lambda_s^{t+u}[X_{si},X_{sj}] \\ \\
&=\displaystyle{\sum_{k=1}^n} m_{ij}^k\lambda_1^{t+u}X_{1k}+...+
\displaystyle{\sum_{k=1}^n} m_{ij}^k\lambda_s^{t+u}X_{sk} \\ \\
&=\displaystyle{\sum_{k=1}^n} m_{ij}^kb_0(X_{1k}+...+X_{sk})+
\displaystyle{\sum_{k=1}^n}
m_{ij}^kb_1(\lambda_1X_{1k}+...+\lambda_sX_{sk}) \\ \\
&+...+ \displaystyle{\sum_{k=1}^n}
m_{ij}^kb_{s-1}(\lambda_1^{s-1}X_{1k}+...+\lambda_s^{s-1}X_{sk}),
\end{array}
$$
showing that $\beta$ is also a $\ZZ$-basis of $\tilde{\ngo}$. Thus
the linear combinations over $\QQ$ of $\beta$ determine a rational
form of $\tilde{\ngo}$, denoted by $\ngoq_{\beta}$, which will be
now showed to be Anosov.  Indeed, it is easy to see that, written in
terms of $\beta$, the hyperbolic automorphism $A$ of $\tilde{\ngo}$
has the form
$$
[A]_{\beta}=\left[\begin{array}{ccc}
B'&&\\
&\ddots&\\
&&B'
\end{array}\right]\in\Gl_{ns}(\ZZ),
$$
where
$$
B'= \left[\begin{array}{ccccc}
0&0&&&-a_0 \\
1&0&&&-a_1\\
0&1&&&\\
&&\ddots&&\\
0&0&&1&-a_{s-1}
\end{array}\right]\in\Gl_s(\ZZ)
$$
is the rational form of the matrix $B$, concluding the proof of the
theorem.
\end{proof}

Different choices of matrices $B$ can eventually give non-isomorphic
Anosov rational forms of $\tilde{\ngo}$, as in the case
$\tilde{\ngo}=\hg_3\oplus\hg_3$ and
$\tilde{\ngo}=\lgo_4\oplus\lgo_4$ (see also \cite{Smn}).  Recall
that two-step nilpotent Lie algebras are graded, so Theorem
\ref{gradedsum} shows that a reasonable classification of Anosov Lie
algebras up to isomorphism is far from being feasible, not only in
the rational case but even in the real case (see \cite{anosov} for
further information).

\begin{remark}\label{exaAnosov} {\rm The only explicit examples of real Anosov Lie
algebras in the literature so far which are not covered by Theorem
\ref{gradedsum} are the following: the free $k$-step nilpotent Lie
algebras on $n$ generators with $k<n$ (see \cite{Dn}, and also
\cite{DkmMlf, Dkm} for a different approach); certain $k$-step
nilpotent Lie algebras of dimension
$d+\binom{d}{2}+...+\binom{d}{k}$ with $d\geq k^2$ (see \cite{Frd});
the $2$-step nilpotent Lie algebra of type $(d,\binom{d}{2}-1)$ with
center of codimension $d$ for $d\geq 5$ (see \cite{DkmDsc}); and the
Lie algebra $\ggo$ (see \cite{anosov}).  Thus $\hg$ is the only new
example over $\RR$ obtained in the classification in dimension $\leq
8$ (see Table \ref{Anosovtable}).  For the known examples of
infranilmanifolds which are not nilmanifolds and admit Anosov
automorphisms we refer to \cite{Shb,Prt,Mlf2}. }
\end{remark}

The {\it signature} of an Anosov diffeomorphism is the pair of
natural numbers $\{ p,q\}=\{\dim{E^+},\dim{E^-}\}$.  It is known
that signature $\{ 1,n-1\}$ is only possible for torus (and their
finitely covered spaces: compact flat manifolds) (see \cite{Frn}).

If $\dim{\ngoq}=n$ then the signature of the Anosov automorphism of
$\tilde{\ngo}^{\QQ}\otimes\RR$
($\tilde{\ngo}^{\QQ}=\ngoq\oplus...\oplus\ngoq$, $s$ times) in the
proof of Theorem \ref{gradedsum} is $\{ np',nq'\}$, $p'+q'=s$, where
$p',q'$ are the numbers of eigenvalues of $B\in\Gl_s(\ZZ)$ having
module greater and smaller than $1$, respectively.  In the
nonabelian case $n$ is necessarily $\geq 3$ and so the signature $\{
2,q\}$ is not allowed for this construction.  We do not actually of
any nonabelian example of signature $\{ 2,q\}$.  We may choose $\{
p',q'\}=\{ 1,s-1\}$ and $\ngoq\otimes\RR=\hg_3$ in order to obtain
signature $\{ 3,3(s-1)\}$ for any $s\geq 2$.

\section{Classification of real Anosov Lie algebras}\label{cla}

We will find in this section all the real Anosov Lie algebras of
dimension $\leq 8$. Our start point is Proposition \ref{coroauto},
which implies that a nonabelian one has to be of dimension $\geq 6$
and gives only a few possibilities for the types in each dimension
$6$, $7$ and $8$.

We use Proposition \ref{auto} to make a few observations on the
eigenvalues of an Anosov automorphism, which are necessarily
algebraic integers.  An overview on several basic properties of
algebraic numbers is given in the Appendix.

\begin{lemma}\label{util}
Let $\ngo$ be a real nilpotent Lie algebra which is Anosov, and let
$A$ and $\ngo=\ngo_1\oplus\ngo_2\oplus\dots \oplus\ngo_r$  be as in
{\rm Proposition \ref{auto}}.  If $A_i=A|_{\ngo_i}$ then the
corresponding eigenvalues $\l_1,\dots,\l_{n_i}$, are algebraic units
such that $1< \degr \l_i \le n_i$ and $\l_1...\l_{n_i}=1$.
\end{lemma}

This follows from the fact that $[A_i]_{\beta_i}\in\Sl_{n_i}(\ZZ)$
and so its characteristic polynomial $p_{A_i}(x) \;\in\; \ZZ[x]$ is
a monic polynomial with constant coefficient $a_0=(-1)^n \det A_i =
\pm 1$, satisfying $p_{A_i}(\l_j)=0$ for all $j=1,\dots,n_i$.

Concerning the degree, it is clear that $\degr \l_j \le n_i$ for all
$j$ and if $\degr \l_j =1$ then $\l_j \in \QQ,$ is a positive unit
and therefore $\l_j=1$, contradicting the fact that $A_i$ is
hyperbolic.

In the following, $\ngo$, $A,$ $A_i$ and $\ngo_i$ will be as in the
previous lemma.  In order to be able of working with eigenvectors,
we will always consider the complex Lie algebra
$\ngo_{\CC}=\ngo\otimes\CC$ and its decomposition
$\ngo_{\CC}=(\ngo_1)_{\CC}\oplus...\oplus(\ngo_r)_{\CC}$, where
$(\ngo_i)_{\CC}=\ngo_i\otimes\CC$.  In the light of Theorem
\ref{abfactor}, we will always assume that $\ngo$ has no abelian
factor.  We now fix more notation that will be used in the rest of
this section.  For simplicity, assume that $\ngo$ is a 2-step
nilpotent Lie algebra.  According to Proposition \ref{auto}, there
exist
$$
\beta_1=\{X_1,X_2,\dots,X_{n_1}\} \qquad {\rm and} \qquad
\beta_2=\{Z_1,Z_2,\dots,Z_{n_2}\},
$$
basis of eigenvectors of $\ngo_1)_{\CC}$ and $(\ngo_2)_{\CC}$ for
$A_1$ and $A_2$, respectively. Let $\l_1,\dots,\l_{n_1}$ and
$\mu_1,\dots,\mu_{n_2}$ be the corresponding eigenvalues.  This
notation will be used throughout all the classification.  The
absence of abelian factor implies that $[\ngo_1,\ngo_1]=\ngo_2$ and
hence we may assume that for each $Z_i$ there exist $X_j$ and $X_l$
such that $Z_i=[X_j,X_l]$. On the other hand, for each $j,l,$ there
exist scalars $a^{j,l}_k\in\CC$ such that $[X_j,X_l]=\sum a_k^{j,l}
Z_k.$ Since $\{Z_k\}$ are linearly independent, for each $k$ we
obtain
\begin{equation}\label{escalares}
\l_j\l_la_k^{j,l}=\mu_ka_k^{j,l}.
\end{equation}
Hence, if $a_k^{j,l}\ne0$, $\mu_k=\l_j\l_l,$ and therefore, if
$a_k^{j,l}\ne 0 \ne a_{k'}^{j,l},$ $\mu_k=\mu_{k'}$. In particular,
if $n_2=2$, since $\mu_1\ne\mu_{2},$ for each $j,l,$ there exist a
unique $k$ such that $[X_j,X_l]=a_k Z_k.$ If it is so, by
(\ref{escalares}), $\l_j\l_l=\mu_k$. When $n_2=3$ the same property
holds. Indeed, $\mu_i \ne \mu_j$  for all $i \ne j $ since $\degr
\mu_i > 1$ for all $i.$

 We are going to consider
all the possible coefficients $a_k^{j,l}$'s only in the cases when
the classification actually leads to a possible Anosov Lie algebra.

\begin{center} {\bf Dimension $\leq 6$} \end{center}

Anosov Lie algebras of dimension $\leq 6$ has already been
classified in \cite{Mlf} and \cite{CssKnnScv}.  We give an
alternative proof here in order to illustrate our approach.

Proposition \ref{coroauto} gives us the following possibilities for
the types of a real Anosov Lie algebra without an abelian factor:
$(3,3)$ and $(4,2)$.

\no {\bf Case $(3,3)$}.  The only real (resp. rational) Lie algebra
of type $(3,3)$ is the free $2$-step nilpotent Lie algebra on $3$
generators $\fg_3$ (resp. $\fg_3^{\QQ}$), which is proved to be
Anosov in \cite{Dn} and  \cite{Dkm,Mlf}.

\no {\bf Case $(4,2)$}.  Let $\ngo$ be a real nilpotent Lie algebra
of type $(4,2)$, admitting a hyperbolic automorphism $A$ as in
Proposition \ref{auto}.  If $\{X_1, \dots X_4\}$ is a basis of
$(\ngo_1)_{\CC}$ of eigenvectors of $A_1$ with corresponding
eigenvalues $\l_1, \hdots, \l_4$, then without any lost of
generality we may assume that we are in one of the following cases:
\begin{enumerate}
\item[(a)] $[X_1,X_2]=Z_1, \qquad [X_1,X_3]=Z_2,$

\item[(b)] $[X_1,X_2]=Z_1, \qquad [X_3,X_4]=Z_2$.

\end{enumerate}

In the first situation, (a) implies that $\l_1^2\l_2\l_3=1,$ and
since $\det A_1=\l_1\l_2\l_3\l_4=1,$ we obtain that $\l_1=\l_4$.
From this it is easy to see that $\degr \l_1=\degr \l_4=2$ and
moreover, $\l_2=\l_3=\l_1^{-1}$ (see Appendix). Therefore, we get to
the contradiction $\mu_1=\mu_2=1$.

Concerning (b), we may assume that there is no more Lie brackets
among the $\{X_i\}$ since otherwise we will be in situation (a), and
thus $\ngo_{\CC} \simeq (\hg_3\oplus\hg_3)_{\CC}$.  This Lie algebra
has two real forms: $\hg_3\oplus\hg_3$ and $\ngoq_{-1}\otimes\RR$
(see paragraph before Proposition \ref{rath3h3}).  The Lie algebra
$\ngoq_{-1}\otimes\RR$ can not be Anosov by Proposition
\ref{region}, (i), and $\hg_3\oplus\hg_3$ is Anosov by Theorem
\ref{gradedsum}.

\begin{center} {\bf Dimension $7$} \end{center}

According to Proposition \ref{coroauto}, if  $\ngo$ is a
$7$-dimensional real Anosov Lie algebra of type
$(n_1,n_2,\dots,n_r)$, then $r=2$ and $\ngo$ is either of type
$(4,3)$ or $(5,2)$. We shall prove that there is no Anosov Lie
algebras of any of these types.

\no {\bf Case $(4,3)$}.  It is easy to see that the eigenvalues of
$A_2$ are three pairs of the form $\lambda_i\lambda_j$, so without
any lost of generality we can assume that two of them are
$\lambda_1\lambda_2$ and $\lambda_1\lambda_3$.  There are four
possibilities for the third eigenvalue of $A_2$, and by using that
$\det{A_1}=1$ and $\det{A_2}= 1$ we get to a contradiction in all
the cases as follows:
\begin{itemize}
\item[(i)] $\l_1\l_2.\l_1\l_3.\l_1\l_4= 1$, then $\l_1^2=1$ contradicting the
hyperbolicity of $A_1$.

\item[(ii)] $\l_1\l_2.\l_1\l_3.\l_2\l_3=1$ implies that $\l_4^2=1$, but then $A_1$
is not hyperbolic.

\item[(iii)] $\l_1\l_2.\l_1\l_3.\l_2\l_4=1$, then $\l_1\l_2= 1$ and so $A_2$ would
not be hyperbolic.

\item[(iv)] $\l_1\l_2.\l_1\l_3.\l_3\l_4=1$, so $\lambda_1\lambda_3=1$ contradicting
the hyperbolicity of $A_2$.
\end{itemize}

\no {\bf Case $(5,2)$}.  Let $\ngo$ be a real nilpotent Lie algebra
of type $(5,2)$, admitting a hyperbolic automorphism $A$ as in
Proposition \ref{auto}.  If $\l_1, \hdots, \l_5,$  are the
eigenvalues of $A_1$ we can either have
\begin{itemize}
\item[(i)] $\l_i \neq \l_j,$ $1 \leq i,j \leq 5$, or \item[(ii)] after reordering if
necessary, $\l_1=\l_2$.
\end{itemize}
Note that in (ii), $\l_1=\l_2$ implies that $ 2\leq 2 \degr \l_1
\leq 5$ and therefore $ \degr \l_1= \degr\l_2 = 2$. From this it is
easy to see that there exist $i \in \{ 3,4,5\}$ such that $\degr
\l_i =1$, contradicting the hyperbolicity of $A_1$. Therefore, we
assume (i).

On the other hand, since $\dim \ngo_2 = 2,$ we have two linearly
independent Lie brackets among the $\{X_i\}$, the basis of
$(\ngo_1)_{\CC}$ of eigenvectors of $A_1$.  Note that if they come
from disjoint pairs of $X_i$, since $\l_1\l_2\l_3\l_4\l_5= 1$, it is
clear that we would have $\l_i= 1$ for some $1 \leq i \leq 5$.
Therefore, without any lost of generality we can only consider the
case when  we have at least the following non trivial Lie brackets:
\begin{equation}\label{eq521}
[X_1,X_2] = Z_1, \qquad\qquad[X_1,X_3] = Z_2.
\end{equation}

In the following we will show that either $X_4$ or $X_5$ are in the
center of $\ngo$, which would generate an abelian factor and hence a
contradiction.  From (\ref{eq521}) we have that
\begin{equation}\label{eq522}
\begin{array}{lcl}
\l_1^2\l_2\l_3= 1, & \text{and then}& \l_4 \l_5 = \l_1.
\end{array}
\end{equation}

Therefore, $[X_4,X_5]=0$ because both of the assumptions
$[X_4,X_5]=cZ_1$ and $[X_4,X_5]=cZ_2$ with $c\ne 0$ leads to the
contradictions $\l_2=1$ and $\l_3=1$, respectively. Also, if
$[X_4,X_j]\ne 0$ and $[X_5,X_k]\ne 0$ for some $1 \leq j,k \leq 3$,
it follows from (i) that we  only have the following possibilities:
$$
\begin{array}{l}
 [X_4,X_3]=c Z_1\qquad  \mbox{and}\qquad  [X_5,X_2]=d Z_2, \;\;\; \mbox{or} \\

 [X_5,X_3]=c Z_1\qquad  \mbox{and}\qquad  [X_4,X_2]=d Z_2,
 \end{array}
 $$
($c,d\ne 0$) which are clearly equivalent. Let us suppose then the
first one, and hence
$$
\begin{array}{lcl}
\text{I.}\; \l_3\l_4=\l_1\l_2 & \text{and} & \text{II.}\;
\l_5\l_2=\l_1\l_3.
\end{array}
$$

 From I, and using that $\l_4\l_5=\l_1$ we obtain $\l_3= \l_2\l_5$.  Therefore by II,
 $\l_1= 1$ which is a
contradiction and then $[X_4,X_j]= 0$ for all $j$ or $[X_5,X_k]= 0$
for all $k$ as we wanted to show.

\begin{center} {\bf Dimension 8} \end{center}

In this case, Proposition \ref{coroauto} gives us the following
possibilities for the types of a real Anosov Lie algebra without an
abelian factor: $(4,4),$ $(5,3),$ $(6,2),$ $(3,3,2)$ and $(4,2,2)$.
Among all this Lie algebras we will show that there is, up to
isomorphism, only three which are Anosov.  One is of type $(4,2,2)$,
one of type $(6,2)$ and one of type $(4,4)$. The first one is an
example of the construction given in \cite{anosov} and Theorem
\ref{gradedsum}, and the second one is isomorphic to \cite[Example
3.3]{anosov}. The last one is a new example.

\no {\bf Case $(4,4)$}.  We will show that there is only one real
Anosov Lie algebra of this type. We first note that there is only
$\binom42=6$ possible linearly independent brackets among the
$\{X_i\}$ and since $\dim [\ngo,\ngo] = 4,$ at most two of them can
be zero. Therefore, without any lost of generality, we can just
consider the following two cases:
\begin{equation}\label{eq441}
 [X_1,X_3]=Z_1, \qquad [X_2,X_4]=Z_2, \qquad [X_2,X_3]=Z_3, \qquad [X_1,X_4]=Z_4,
\end{equation}
that is, the possible zero brackets corresponds to disjoint pairs of
$\{X_i\}$ (namely $\{X_1,X_2\}$ and $\{X_3,X_4\})$; and the other
case is
\begin{equation}\label{eq442}
 [X_1,X_4]=Z_1, \qquad [X_2,X_4]=Z_2, \qquad [X_3,X_4]=Z_3, \qquad [X_2,X_3]=Z_4,
\end{equation}
corresponding to the case of non disjoint pairs, $\{X_1,X_2\}$ and
$\{X_1,X_3\}.$

However, the second case is not possible because we would have
$$
\text{I)}\; \l_1\l_2\l_3\l_4=1 \qquad \text{and} \qquad \text{II)}
\; \l_1\l_2^2\l_3^2\l_4^3=1.
$$
It follows that $\l_2\l_3= \l_4^{-2}$ and replacing this in I) we
get $\l_1=\l_4$. This implies that the $\l_i$'s have all degree two,
and $\l_2 = \l_3 = \l_4^{-1}$ (see Appendix). Hence $\mu_3=\l_3\l_4=
1$, contradicting the hyperbolicity of $A_2$.

Concerning case (\ref{eq441}), if we assume that $[X_1,X_2]=0$ and
$[X_3,X_4]=0$
$$
\begin{array}{lclcl}
A= \left[ \begin{smallmatrix}A_1& \\
               &A_2 \end{smallmatrix}\right], & \text{where}&

A_1= \left[ {\begin{smallmatrix}\l &      &    & \\
                                  &\l^{-1}&    & \\
                                  &       &\l^2&  \\
                                 &       &     &\l^{-2}\end{smallmatrix}}\right]
&\text{and}&
A_2= \left[ \begin{smallmatrix}\l^3&      &    & \\
                                  &\l^{-3}&    & \\
                                  &       &\l&  \\
                                 &       &     &\l^{-1}\end{smallmatrix}\right]
\end{array}
$$

is an automorphism of $\ngo$ for any $\l \in \RR^*$.  If
$\l\in\RR^*$ is an algebraic integer such that $\l +\l^{-1} = 2a$,
$a \in \ZZ$, $a\geq 2$, then it is easy to check that
\begin{equation}\label{base}
\begin{array}{ll}
\beta=&\left\{ X_1+X_2,\; (a^2-1)^{\unm}(X_1-X_2),\; X_3+X_4,\; (a^2-1)^{\unm}(X_3-X_4),\right.\\
 &\left.Z_1+Z_2,\;
(a^2-1)^{\unm}(Z_1-Z_2),\; Z_3+Z_4,\;
(a^2-1)^{\unm}(Z_3-Z_4)\right\}
\end{array}
\end{equation}
is a $\ZZ$-basis of $\ngo$. Moreover, if $B=
\left[\begin{smallmatrix}
a&a^2-1\\
1&a
\end{smallmatrix}\right]$, then the matrix of $A$ in terms of the basis $\beta$ is given by
$$
\left[A\right]_{\beta}= \left[\begin{array}{cccc}
B &&& \\
& B^{2}&&\\
&&B^3&\\
&&&B
\end{array}\right]\in\Sl(8,\ZZ),
$$
showing that $\ngo$ is Anosov.  Recall that this $\ngo$ is
isomorphic to the Lie algebra $\hg$ given in Example \ref{dual42}.

It follows from Scheuneman duality that there is only one more real
form of $\hg_{\CC}$, namely, the dual of the Lie algebra
$\ngoq_{-1}\otimes\RR$ of type $(4,2)$ ($\hg$ is dual of
$\hg_3\oplus\hg_3$).  The fact that such a Lie algebra is not Anosov
will be proved in Section \ref{cla2}, Case $\hg$.

We will now show that if we add any more nonzero brackets in case
(\ref{eq441}), then the new Lie algebra $\tilde{\ngo}$ does not
admit a hyperbolic automorphism any longer.  Suppose then that
 $$0 \ne [X_1,X_2]=a_1Z_1 + a_2Z_2 + a_3Z_3 + a_4Z_4.$$
As we have already pointed out at the beginning of this
classification, since A is an automorphism and $Z_i$ are linearly
independent, it follows that if $a_j \ne 0$ then $\l_1\l_2=\mu_j.$
Therefore, at most two of them can be non zero.

If $[X_1,X_2]=a_jZ_j$ then we can change $Z_j$ by
$\tilde{Z_j}=a_jZ_j$ and the corresponding bracket in (\ref{eq441})
by $[X_1,X_2]$ and we will be in the conditions of case
(\ref{eq442}).

If $[X_1,X_2]=a_jZ_j+a_kZ_k,$  $a_j, \; a_k \ne 0,$ then we have
that $\l_1\l_2=\mu_j=\mu_k.$ One can check that for all the choices
of $j,k$ we obtain $\l_i=\l_r=\l_s$ for some $1 \le i,r,s \le 4$
which is not possible because it implies that $2 \le 3\degr \l_i \le
4$ and then $\degr \l_i=1.$

Hence we get $[X_1,X_2]=0$ and by using the same argument we also
obtain $[X_3,X_4]=0$ as we wanted to show.

We also note that for any choice of nonzero scalars $a,b,c,d$, the
Lie algebra $\tilde{\ngo}$ given by
$$
 [X_1,X_3]=aZ_1 \qquad [X_2,X_4]=bZ_2 \qquad [X_2,X_3]=cZ_3 \qquad [X_1,X_4]=dZ_4,
$$
is isomorphic to $\ngo.$

\no {\bf Case $(5,3)$}.  We shall prove that there are no Lie
algebras of this type with no abelian factor
 admitting a hyperbolic automorphism.

Suppose that $A$ is as in Proposition \ref{auto}. Hence as we have
already pointed out, the eigenvalues of $A_1$, $\l_1,\dots,\l_5$ are
algebraic integers with $2\le \degr \l_j \le 5$ for all $1\le j\le
5$. As we have seen in case (5,2), we can assume that $\l_i \ne
\l_j$ for all $i\ne j$ since otherwise we will have
 that there exists $ k$ with $\degr \l_k =1$, contradicting the hyperbolicity of
$A_1$.  In this situation it is easy to see that
\begin{equation}\label{eq531}
\text{if}\; \sharp\left(\{X_i,X_j\} \cap \{X_k,X_l\}\right)=1 \qquad
\text{ then}\qquad [X_i,X_j] \notin \CC [X_k,X_l].
\end{equation}
Moreover, since $2\le \degr \mu_k \le 3$ we have that $\mu_k \ne
\mu_l$ for all $1 \le k\ne l \le 3$ and then for all $i,j$ there
exist $k$ such that $[X_i,X_j] \in \CC Z_k.$

On the other hand, it is clear that we can split the set of Lie
algebras of this type according to the following condition:
\begin{equation}\label{cond53}
\begin{array}{l}
\text{There are two disjoint pairs of}\;\; \{X_i\} \;\; \text{such that the corresponding } \\
\text{Lie brackets are linearly independent.}
\end{array}
\end{equation}

Note that if $\ngo$ does not satisfy this condition, we will have
that

\begin{equation}\label{obs53}
\{X_i,X_j\}\cap \{X_l,X_k\} = \emptyset \quad \Rightarrow \quad
[X_i,X_j] \in \CC [X_l,X_k].
\end{equation}

If (\ref{obs53}) holds, we can assume without any lost of generality
that
\begin{equation}\label{eq532}
\begin{array}{ll}
[X_1,X_2]=Z_1 & [X_1,X_3]=Z_2,
\end{array}
\end{equation}
and for $Z_3$ we have two possibilities
$$
{\rm a)}\; [X_1,X_4]=Z_3, \qquad {\rm b)}\;[X_2,X_3]=Z_3.
$$
We will now show that any of this assumptions leads to a
contradiction.

Concerning a), we have that $[X_5,X_k] \ne 0 $ for some $1\le k\le
4,$ but since $\l_i\ne\l_j,$ when $i\ne j$ it is clear that $k\ne1.$
We can assume then that $k=2$, since every other choice (i.e.
$k=3,4$) is entirely analogous. Now, since $\{5,2\} \cap \{1,3\} =
\emptyset,$ by (\ref{obs53}) we have that
 $[X_5,X_2]\in \CC Z_2$, and analogously,
$\{5,2\} \cap \{1,4\} = \emptyset$ and then $[X_5,X_2] \in \CC Z_3$,
giving the contradiction $[X_5,X_2]=0$.

In case b) $\l_1\l_2\l_3=1$, and therefore $\l_4\l_5=1$.  Thus
$[X_5,X_4] = 0$, and we may assume that $0\ne [X_4,X_1]\in\CC Z_3$
and $0\ne [X_5,X_2]\in\CC Z_2$. Therefore,   $\l_5\l_2=\l_1\l_3$ and
$\l_4\l_1=\l_2\l_3$, and since $\l_4\l_5=1$, we get to the
contradiction $\l_3=1$.

We can assume then that $\ngo$ satisfies condition (\ref{cond53})
and thus without any lost of generality we can suppose that
 \begin{equation}\label{eq533}
 [X_1,X_2]=Z_1 \qquad [X_3,X_4]=Z_2.
\end{equation}
Note that we can not have $[X_5,X_j]=Z_3$ because this would imply
$\l_j=1$ by using that $\l_1\dots\l_5=1$. Let us say then that
$[X_5,X_j]=aZ_1$, $a\ne 0$.  From (\ref{eq531}) we have that $j \ne
1,2,$ and since both cases $j=3$ and $j=4$ are completely analogous,
we will just analyze the case $j=3$. This is
$$
[X_1,X_2]=Z_1,\;\; [X_3,X_4]=Z_2,\;\; [X_5,X_3]=aZ_1.
$$
Also, since $Z_3 \in  [\ngo,\ngo]$ there is $1\le k,k' \le 4$ such
that $[X_k,X_{k'}]=Z_3$, and by the above observations, it is easy
to see that
$$\begin{array}{ll}
\{k,k'\}=&\left\{ \begin{array}{l} \{1,3\} \text{ or (equivalently) }  \{2,3\} \\ \\
                                     \{1,4\} \text{ or (equivalently) }  \{2,4\}
                 \end{array} \right.
                 \end{array}$$

To finish the proof, we will see that both cases leads to a
contradiction.  The idea is to show that one of the $\l_i$ is equal
to one of the $\mu_j,$ and since the conjugated numbers are uniquely
determined, this implies that every $\mu_j$ appears as a $\l_k.$
From here it is easy to check in both cases that this is not
possible.

Indeed, if $[X_1,X_3]=Z_3$, since $ 1 = \l_5\l_3\l_3\l_4\l_1\l_3,$
we have that $\l_3^2=\l_2$. Therefore, $\l_5\l_3=\l_1\l_2=
\l_1\l_3^2$ and so $\l_5 = \l_1\l_3 = \mu_3$.  Hence, there exists
$i$ such that $\mu_1=\l_1\l_3^2= \l_i$.  It is clear that $i \ne
1,2,3,5$ and if $\l_1\l_3^2=\l_4,$ since
$1=\l_1\l_2\l_3\l_4\l_1\l_3= \l_1^2\l_3^4\l_4,$ then $ 1 =
\l_1^3\l_3^6= \mu_1^3$ contradicting the fact that $A_2$ is
hyperbolic.

\vspace{.3cm}

Now, if $[X_1,X_4]=Z_3$, then
\begin{itemize}
\item[(i)] $ 1 = \l_5\l_3\l_3\l_4\l_1\l_4,$ and from there $\l_2 = \l_3\l_4 =
\mu_2,$ and

\item[(ii)]  $1 = \l_1\l_2\l_3\l_4\l_1\l_4,$ hence $\l_5 = \l_1\l_4 = \mu_3.$
\end{itemize}

Therefore, as we have observed before, there is $k$ such that
$\mu_1=\l_k.$ This implies that $\l_1\l_2 = \l_5\l_3 = \l_k$ for
some $1\le k\le 5$. Again, it is clear that $k \ne 1,2,3,5,$ and if
$\l_1\l_2=\l_5\l_3= \l_4,$ then by (ii) $ \l_1\l_4\l_3= \l_5\l_3 =
\l_4$  and hence $\l_1\l_3= 1.$ From this, using that $1=\det
A_2=\l_4\l_2\l_5$, we obtain that $\l_1\l_2=\l_5\l_3=
\frac{1}{\l_2\l_4}.\frac{1}{\l_1},$ or equivalently
$\l_4^2=(\l_1\l_2)^2= \frac{1}{\l_4}$ and then $\l_4= 1$
contradicting the fact that $A_1$ is hyperbolic, and concluding the
proof of case $(5,3)$.

\no {\bf Case $(6,2)$.}  We will prove in this case that there is,
up to isomorphism, only one Anosov Lie algebra with no abelian
factor. As usual, let $A$ be an Anosov automorphism of $\ngo$ and
$\{X_1,\dots,X_6,Z_1,Z_2\}$ a basis of $\ngo_{\CC}$ of eigenvectors
of $A$, $\l_1,\dots,\l_6,\mu_1,\mu_2$ the eigenvalues as above.

As we have mentioned before, since $\mu_1 \ne \mu_2,$ for all $i,j$
there exists $k$ such that $[X_i,X_j]\in\CC Z_k$.  Also, if
$\dim{[X_i,(\ngo_1)_{\CC}]}=1$ for any $i$, then $\ngo_{\CC}$ is
either isomorphic to $(\hg_3\oplus\RR\oplus\hg_3\oplus\RR)_{\CC}$ or
$(\hg_3\oplus\hg_5)_{\CC}$.  The first one has two real forms:
$\hg_3\oplus\hg_3\oplus\RR^2$ and
$(\ngoq_{-1}\otimes\RR)\oplus\RR^2$, of which only
$\hg_3\oplus\hg_3\oplus\RR^2$ is Anosov by Theorem \ref{abfactor}
and the classification in dimension $6$.  The only real form of
$(\hg_3\oplus\hg_5)_{\CC}$ is $\hg_3\oplus\hg_5$ (see Remark
\ref{frh3h5}), and by Proposition \ref{rath3h5}, $\hg_3\oplus\hg_5$
has only one rational form with Pfaffian form $f(x,y)=xy^2$.  It
then follows from Proposition \ref{region}, (ii) that it is not
Anosov.

Therefore, we can assume that
\begin{equation}\label{eq621}
[X_1,X_2]=Z_1, \qquad [X_1,X_3]=Z_2.
\end{equation}
From this, one has that
\begin{equation}\label{eq622}
\l_1^2\l_2\l_3= 1,\quad \mbox{or equivalently} \quad
\l_1=\l_4\l_5\l_6.
\end{equation}

In what follows, we will first show that there exist a reordering
$\beta$ of $\{X_1,...,X_6\}$, such that
\begin{equation}\label{assertion}
[A_1]_\beta=\left[\begin{smallmatrix}\l&&&&& \\&\l^{-1}&&&&\\&&\nu&&&\\&&&\nu^{-1}&&\\
&&&&\mu&\\
&&&&&\mu^{-1}\end{smallmatrix}\right],
\end{equation}
and after that we will see that this implies that $\ngo_{\CC} \simeq
\ggo_{\CC},$ the complexification of the Lie algebra defined in
(\ref{62alg}), which is proved to be Anosov in \cite[Example
3.3]{anosov}.  Moreover, $\ggo$ is known to be the only real form of
$\ggo_{\CC}$ (see Remark \ref{frg}).

To do this, let us first assume that
\begin{itemize}
\item[a)] $\l_i=\l_l,$ denoted by $ \l,$ for some $1 \le i \ne j \le 6$.
\end{itemize}

Thus $\degr \l = 2,$ or $\degr \l = 3$, but $\degr \l=3$ is not
possible. In fact, if $\degr \l=3$ then there exist a reordering of
$\{ X_i\}$ such that the matrix of $A_1$ in the new basis is
$$
A_1= \left[\begin{array}{cc} B&0\\ 0&B
\end{array}\right],\quad \mbox{ where} \quad B=\left[\begin{smallmatrix} \l &&\\ &
\mu&\\&&({\l\mu})^{-1}\end{smallmatrix}\right]
$$
is conjugated to an element in $\Sl_3(\ZZ)$.  This says that
$\l_1,\l_2,\l_3 \in \left\{\l, \mu,(\mu\l)^{-1}\right\},$ and using
(\ref{eq622}) one can see that $\l_1= \l_2$ (or equivalently
$\l_1=\l_3$), since every other choice ends up in a contradiction.
Therefore, we may assume that $\l_1=\l_2=\l$ and so $\l_3=
\l^{-3}=\mu$.  Since every eigenvalue of $A_1$ has multiplicity $2$,
we have that after a reordering if necessary, $\l_4= \l_5=\l^2$ and
$\l_6=\l^{-3}$. Therefore, the matrix of $A$ in the basis $\beta=\{
X_1,X_2,...,X_6,Z_1,Z_2\}$ is given by
$ [A]_\beta=\left[\begin{smallmatrix}A_1& \\
&A_2\end{smallmatrix}\right]$ where
$$
A_1=\left[\begin{smallmatrix}\l&&&&& \\&\l&&&&\\&&\l^{-3}&&&\\&&&\l^2&&\\&&&&\l^2&\\
&&&&&\l^{-3}\end{smallmatrix}\right]  \qquad \text{and} \qquad
 A_2=\left[\begin{smallmatrix}\l^2& \\&\l^{-2}
\end{smallmatrix}\right].
$$
Hence, since $A$ is an automorphism of $\ngo$, one gets that
$X_4,X_5 \in(\zg \cap \ngo_1)_{\CC}$, contradicting our assumption
of no abelian factor.  Thus $\degr \l=2,$ from where assertion
(\ref{assertion}) easily follows.

On the other hand, if
\begin{itemize}
\item[b)] $\l_i\ne \l_j$ for all $ i \ne j,$
\end{itemize}
with no loss of generality, we can assume that $[X_4,X_j] =aZ_1$,
$a\ne 0$, for some $j\in\{ 3,5,6\}$.

If $j=5$ then it follows from $1=\det A_2= \l_4\l_5\l_1\l_3,$ that
\begin{equation}\label{eq623}
\l_2\l_6=1.
\end{equation}
Now, we also have that $[X_6,X_k] \ne 0$ for some $k$, and hence it
is easy to see that we can either have
\begin{enumerate}
\item[I)] $[X_6,X_3] = bZ_1$, $b\ne 0$, or

\item[II)] $[X_6,X_4] =cZ_2$, $c\ne 0$, (or equivalently $[X_6,X_5]=cZ_2$).
\end{enumerate}
In case I), $\l_6\l_3=\l_1\l_2$ and so by (\ref{eq623}) we have that
$\l_3=\l_1\l_2^2$.  By using (\ref{eq622}) we get to the
contradiction $\mu_1=1.$

Concerning II), since $\l_6\l_4=\l_1\l_3,$ we obtain from
(\ref{eq622}) that $\l_5\l_3=1$ and therefore $\l_1\l_4=1$.  This
together with (\ref{eq623}) implies assertion (\ref{assertion}). The
case when $j=6$ is entirely analogous to the case $j=5$ and so we
are not going to consider it.

If $j=3$ then $\l_4\l_3=\l_1\l_2$ and by (\ref{eq622}) it is easy to
see that
\begin{equation}\label{eq624}
\l_3=\l_5\l_6\l_2.
\end{equation}
Analogously to the previous case, since $[X_5,X_k] \ne 0$ for some
$k$, it is easy to see that we can either have
\begin{enumerate}
\item[I)] $[X_5,X_6] =b\,Z_1,$ or

\item[II)] $[X_5,X_2] =c\,Z_2$ (or equivalently $[X_5,X_4]=cZ_2$).
\end{enumerate}
It is easy to deduce from the situation I) that (\ref{eq624})
implies that $\mu_2=\mu_1^2$ and so both of them are equal to $1$
contradicting the fact that $A_2$ is hyperbolic.

In case II), $\l_5\l_2=\l_1\l_3$ and it follows from (\ref{eq624})
that $\l_6\l_1=1$.  Also, since $\ngo$ has no abelian factor, it is
easy to see that $[X_6,X_4] = dZ_2$, $d\ne 0$, and therefore
$\l_6\l_4=\l_1\l_3.$ Hence, using (\ref{eq622}) we obtain
$\l_2\l_4=1,$ from where assertion (\ref{assertion}) follows.

\vs

To finish the proof we must study the case when (\ref{assertion})
holds, that is,
$$
A_1= \left[\begin{smallmatrix} A_{\l}&&\\&A_{\nu}&\\&&A_\mu
\end{smallmatrix}\right]\;\; \text{where} \;\;A_\eta=\left[\begin{smallmatrix}\eta &\\ &
\eta^{-1}\end{smallmatrix}\right].
$$
Let $\l_1 =\l,$ $ \l_2=\nu$ and thus, by (\ref{eq622}),
$\l_3=\frac{1}{\l^2\nu}.$ It is easy to see that $\l_3$ is different
from $\l^{-1}$ or $\nu^{-1}.$ Therefore, after a reordering if
necessary, we have that
$$
A_1= \left[\begin{smallmatrix}
\l&&&&&\\&\nu&&&&\\&&(\l^2\nu)^{-1}&&&\\&&&\l^2\nu&&\\&&&&\l^{-1}&\\&&&&&\nu^{-1}
\end{smallmatrix}\right]\;\; \text{and} \;\;A_2=\left[\begin{smallmatrix}\l\nu &\\ &
(\l\nu)^{-1}\end{smallmatrix}\right].
$$
Using that $A$ is an automorphism, one can see that $[V_1,V_2]=0$,
where $V_1=\la X_1,X_2,X_3\ra_{\CC}$ and $V_2=\la
X_4,X_5,X_6\ra_{\CC}$.  Moreover, since $\ngo_{\CC}$ has no abelian
factor $[V_1,V_2]=\la Z_1,Z_2\ra_{\CC}$.  From the classification of
$2$-step nilpotent Lie algebras with $2$-dimensional derived algebra
in terms of Pfaffian forms given in Section \ref{rational}, it
follows that there is only one Lie algebra satisfying these
conditions and so $\ngo_{\CC}$ is isomorphic to $\ggo_{\CC}$, as was
to be shown.

\no {\bf Case $(3,3,2)$}.  We will show in this case that there is
no Anosov Lie algebra. We will begin by noting that since $\ngo_2$
has dimension three, we may assume that
$$
[X_1,X_2]=Y_3, \quad [X_1,X_3]=Y_2, \quad [X_2,X_3]=Y_1,
$$
where $\{ X_1,X_2,X_3\}$ and $\{ Y_1,Y_2,Y_3\}$ are basis of
$(\ngo_1)_{\CC}$ and $(\ngo_2)_{\CC}$ of eigenvectors of $A_1$ and
$A_2$, respectively.

It follows that
\begin{equation}\label{eq3321}
[X_1,Y_1]=0, \quad [X_2,Y_2]=0, \quad [X_3,Y_3]=0,
\end{equation}
since any of them would be an eigenvector of $A$ of eigenvalue
$\l_1\l_2\l_3 =1$ and then $A_3$ would not be hyperbolic.

On the other hand, since $Z_1,Z_2 \in \ngo_3$ we have that for some
$i,j,k,l$
$$
[X_i,Y_j]=Z_1, \qquad [X_k,Y_l]=Z_2,
$$
and thus $i\ne k$.  Indeed, if $i=k$ then $j\ne l$ and by
(\ref{eq3321}) $j,l\ne i$. This would imply that
$\l_i.\l_i\l_j.\l_i.\l_i\l_l=1$ and so $\l_i^3=1$, a contradiction.

Hence we can assume that
$$
[X_1,Y_j]=Z_1 \qquad [X_2,Y_l]=Z_2.
$$
For the pairs $(j,l)$ we have four possibilities as follows:
$(2,1),(2,3),(3,1)$ and $(3,3)$. In order to discard some of them,
we recall that since $\dim \ngo_1=3,$ $\l_i\ne \l_j$ for all $1\le
i,j \le 3$ and from this, it follows that $(j,l)\ne(3,1)$ or
$(2,3).$ Indeed, if $(j,l)=(3,1)$ (or $(2,3)$) we have that
$\l_1\l_1\l_2 \l_2\l_2\l_3 =1$ (or $\l_1\l_1\l_3\l_2\l_1\l_2=1$).
Hence $ \l_1\l_2^2=1$ (or $\l_1^2\l_2=1$) and we get to the
contradiction $\l_2=\l_3$ (or $\l_1=\l_3$).

It is also easy to see that $(j,l)\ne (3,3)$ since this implies
$\l_1\l_1\l_2\l_2\l_1\l_2$ and so $\l_1\l_2=1$, contradicting the
fact that $A_2$ is hyperbolic.  Finally, assume that $(j,l)=(2,1),$
that is, in $\ngo_{\CC}$ we have at least the following non trivial
brackets:
\begin{equation}\label{eq3322}
\begin{array}{lcl}
[X_1,X_2]=Y_3,&[X_1,X_3]=Y_2,&[X_2,X_3]=Y_1,\\&&\\

[X_1,Y_2]=Z_1,&[X_2,Y_1]= Z_2.&
\end{array}
\end{equation}
Let $\l_1=\l$ and $\l_2=\nu,$ then the matrix of $A$ is given by
$$
\begin{array}{lcl}
[A]=\left[\begin{smallmatrix}
B&&&\\&B^{-1}&&\\&&\frac{\l}{\nu}&\\&&&\frac{\nu}{\l}\end{smallmatrix}\right],&\quad
\text{where}\quad&
B=\left[\begin{smallmatrix}\l&&\\&\nu&\\&&\frac{1}{\l\nu}\end{smallmatrix}\right],
\end{array}
$$
and $B$ is conjugated to an element of $\Sl_3(\ZZ)$. Thus
$\frac{\l}{\nu}$ is an algebraic unit with $|\frac{\l}{\nu}|\ne 1$
and  $\degr{\frac{\l}{\nu}}=2$.  It is easy to see that under such
conditions $\frac{\l}{\nu}$ is necessarily a real number.  Since the
possibilities for $\nu$ are either $\nu=\overline{\l}$ or
$\frac{1}{|\l|^2}$, we obtain that $\l,\nu\in\RR$, which is a
contradiction by the following lemma applied to $\l^2,\nu^2$.
 This concludes the proof of this case.

\begin{lemma} {\rm Let $\l_1 ,\l_2$ be two positive totally real algebraic integers of degree 3. If $\l_1$ and $\l_2$
are conjugated and units then $\frac{\l_1}{\l_2}$ can never have
degree two.}
\end{lemma}
\begin{proof} Let $\l_1$ and $\l_2$ be as in the lemma, then the minimal polynomial of $\l_i$ is given by
$m_{\l_i}(x)=(x-\l_1)(x-\l_2)(x-\l_3),$ where $\l_1\l_2\l_3=\pm 1$.
Since $m_{\l_i}$ has its coefficients in $\ZZ,$ we have that
$$
\l_1+\l_2+\l_3\in \ZZ, \qquad \frac{1}{\l_1}+
\frac{1}{\l_2}+\frac{1}{\l_3}\in \ZZ,
$$
and hence
\begin{equation}\label{eq333}
\l^2_1+\l^2_2+\l^2_3=d \in \ZZ.
\end{equation}
On the other hand, if we assume that $\l_1/\l_2$ has degree two then
$\frac{\l1}{\l_2}+\frac{\l_2}{\l_1}=a\in\ZZ$, and thus
$$
\frac{\l_1}{\l_2}=\frac{a}{2}+\sqrt{\frac{a^2}{4}-1} \quad
\text{and}\quad
\frac{\l_2}{\l_1}=\frac{a}{4}-\sqrt{\frac{a^2}{2}-1}.
$$
Recall that $a\geq 2$.  We also note that
$\frac{\l_1}{\l_2}=\pm\l_1^2\l_3$ and
$\frac{\l_2}{\l_1}=\pm\l_2^2\l_3$, and hence $\l_1^2 = \pm
\frac{1}{\l_3}\left( \frac{a}{2}+\sqrt{\frac{a^2}{4}-1}\right)$ and
$\l_2^2
=\pm\frac{1}{\l_3}\left(\frac{a}{2}-\sqrt{\frac{a^2}{4}-1}\right)$.
By replacing this in (\ref{eq333}) we obtain $\pm\frac{1}{\l_3}a +
\l_3^2 = d$, or equivalently,
$$
 \l_3^3 - \l_3 d \pm a = 0.
$$
This means that $p(x)=x^3-dx\pm a\;$ is a monic polynomial of degree
3 with coefficient in $\ZZ$ which is annihilated by $\l_3$. Hence it
is equal to the minimal polynomial of $\l_3$ and then $a= \pm 1$,
which is a contradiction since as we have observed above, $a\geq 2$.
\end{proof}

We would like to point out that in this lemma, we are strongly using
the fact that $\l_1$ and $\l_2$ are totally real algebraic numbers
and units. Indeed, if we consider $p(x)=x^3-2$, the roots of $p$ are
$\left\{ \l_1=2^{1/3}, \; \l_2=\omega 2^{1/3},\; \l_3= \omega^2
2^{1/3}\right\},$ where $\omega^2+\omega+1=0$. Since $x^3-2$ is
indecomposable over $\QQ$, we have that $\degr \l_i=3$ for all
$i=1,2,3,$ and however $\l_2. \frac{1}{\l_1}=\omega$ has degree two.

\no {\bf Case $(4,2,2)$}.  Let $\ngo$ be a nilpotent Lie algebra of
type $(4,2,2)$ and let $A$ be an hyperbolic automorphism with
eigenvectors $\{X_1, \dots , X_4, Y_1,Y_2,Z_1,Z_2\}$, a basis of
$\ngo_{\CC}$, and corresponding eigenvalues $\l_1,
\dots,\l_4,\eta_1,\eta_2,\mu_1,\mu_2$ as in Proposition
\ref{coroauto}.

 Since $\eta_i = \l_{j}\l_{k}$ we have the
following two possibilities:
\begin{enumerate}
\item[(I)] In the decomposition of $\eta_1\eta_2$ as product of $\l_i$ at least one
of the $\l_i$ appears twice, or

\item[(II)] $\eta_1=\l_1\l_2$, $\eta_2=\l_3\l_4,$ and $\l_i\ne\l_j$ for $1\le i,j
\le 4$.
\end{enumerate}

\vspace{.1cm}

In the first case we can either have
$$\begin{array}{lcl}
\text{(a)}\; \eta_1=\l_1\l_2, \;\; \eta_2=\l_1\l_3, & \text{(b)}\;
\eta_1=\l_1^2, \;\; \eta_2=\l_2\l_3,\;\; \text{or}\; &\text{(c)}\;
\eta_1=\l_1^2, \;\; \eta_2=\l_1^{-2}.
\end{array}$$
Note that (a) and (b) implies that $\l_1\l_2.\l_1\l_3 =1$ and hence
$\l_4=\l_1.$ Thus $\degr \l_4 = \degr \l_1=2,$ and moreover,
$\l_2=\l_3=\pm\l_1^{-1}.$ Therefore in case (a) we get to  the
contradiction $\eta_1=\eta_2=\pm 1$, and case (b) becomes (c).

So it remains to study case (c).  There is no lost of generality in
assuming that $\l_1=\l_2=\l$ and $\l_3=\l_4=\l^{-1}$ and from this,
using the Jacobi identity, it is easy to see that the possible
nonzero brackets are
\begin{equation}\label{eq4220}
\begin{array}{lll}
[X_1,X_2]=Y_1,& [X_2,Y_1]=a\;Z_1&  [X_1,Y_1]=a'\;Z_1 \\ &&\\

[X_3,X_4]=Y_2, &  [X_3,Y_2]=b\;Z_2. & [X_4,Y_2]=b'\;Z_2.
\end{array}
\end{equation}

Since $\ngo_{\CC}$ has no abelian factor, we have that $a\ne 0$ or
$a'\ne 0$ and $b \ne 0$ or $b'\ne 0$. Let $\ngo_0$ be the ideal of
$\ngo_{\CC}$ generated by $\{X_1,X_2,Y_1,Z_1\}$ and $\ngo_0'$ the
ideal generated by $\{X_3,X_4,Y_2,Z_2\}.$ By the above observation,
they are both four dimensional 3-step complex nilpotent Lie
algebras. It is well known that there is up to isomorphism only one
of such Lie algebras and therefore $\ngo_0$ and $\ngo_0'$ are both
isomorphic to $(\lgo_4)_{\CC}$ and
$\ngo_{\CC}=(\lgo_4\oplus\lgo_4)_{\CC}$.  By Remark \ref{frl4l4}, we
know that $\lgo_4\oplus\lgo_4$ is the only real form of
$(\lgo_4\oplus\lgo_4)_{\CC}$, and it is Anosov by Theorem
\ref{gradedsum}. This concludes case (I).

We will now study case (II).  We can assume that
\begin{equation}\label{eq4221}
[X_1,X_2]=Y_1, \qquad  [X_3,X_4]=Y_2.
\end{equation}
Moreover, due to our assumption it is easy to see that there is no
more non-trivial Lie brackets among them.  On the other hand, we
have that $Z_i \in \ngo_3$ and then for each $i=1,2$
$$Z_i=[X_{j_i},Y_{k_i}].$$
If $k_1=k_2$ we may assume that $k_1=k_2=1$.  By using Jacobi
identity and the previous observation, one can see that
$j_1,j_2\notin\{ 3,4\}$, and hence we get
$$
[X_1,Y_1]=Z_1, \qquad [X_2,Y_1]=Z_2.
$$
From this we have that $\l_1.\l_1\l_2.\l_2.\l_1\l_2=1$ and therefore
$\l_1^3\l_2^3=1$, a contradiction.

Otherwise, we can assume that $k_1=1$ and $k_2=2$. Therefore
$\l_{j_1}\l_1.\l_{j_2}\l_2.\l_3\l_4=1$ and then
$\l_{j_1}\l_{j_2}=1.$ Hence $j_1\ne j_2$ and since $\l_1\l_2 \ne 1$
and $\l_3\l_4\ne 1$ we can suppose that $\l_1\l_3=1$ and
$\l_2\l_4=1$.  Without any lost of generality we can assume that
\begin{equation}\label{eq4222}
[X_1,Y_1]=Z_1  \quad \text{and} \quad [X_3,Y_2]=Z_2,
\end{equation}
since by Jacobi $[X_1,Y_2]=[X_3,Y_1]=0$.  Note that we have obtained
that the matrix of $A$ is given by
$$
\begin{array}{llll}
[A_1]=\left[\begin{smallmatrix}\l&&&\\&\nu&&\\&&\l^{-1}&\\&&&\nu^{-1}\end{smallmatrix}
\right],& [A_2]=\left[\begin{smallmatrix}\l\nu& \\ &
(\l\nu)^{-1}\end{smallmatrix}\right]& \text{and}\quad
[A_3]=\left[\begin{smallmatrix} \l^2\nu&\\&(\l^2\nu)^{-1}
\end{smallmatrix}\right].
\end{array}
$$
From this, since $\l\ne \nu$ and $A\in\Aut(\ngo_{\CC})$, it is easy
to see that we can not have other nonzero Lie brackets on
$\ngo_{\CC}$ but (\ref{eq4221}), (\ref{eq4222}), $[X_1,X_4]=aZ_1$
and $[X_2,X_3]=bZ_2$.  This Lie algebra is isomorphic to the one
with $a=b=0$ (by changing for $\tilde{X}_4=X_4+Y_1$,
$\tilde{X}_2=X_2+Y_2$), and then $\ngo_{\CC}$ is again isomorphic to
$(\lgo_4\oplus\lgo_4)_{\CC}$.

\vspace{.5cm}

We summarize the results obtained in this section in the following

\begin{theorem}\label{clasreal}
Up to isomorphism, the real Anosov Lie algebras of dimension $\le 8$
are: $\RR^n,$ $n=2,\dots,8,$ $\hg_3\oplus\hg_3,$ $\fg_3,$
$\hg_3\oplus\hg_3\oplus\RR^2,$ $\fg_3\oplus\RR^2,$ $\ggo,$ $\hg,$
and $\lgo_4\oplus\lgo_4.$
\end{theorem}

\section{Classification of rational Anosov Lie algebras}\label{cla2}

In Section \ref{cla}, we have found all real Lie algebras of
dimension $\leq 8$ having an Anosov rational form (see Theorem
\ref{clasreal} or Table \ref{Anosovtable}).  On the other hand, the
set of all rational forms (up to isomorphism) for each of these
algebras has been determined in Section \ref{rational} (see Table
\ref{ratforms}).  In this section, we shall study which of these
rational Lie algebras are Anosov, obtaining in this way the
classification in the rational case up to dimension $8$.

\no {\bf Case $\fg_3$ (type $(3,3)$)}.  There is only one rational
form $\fg_3^{\QQ}$ in this case which is proved to be Anosov in
\cite{Dn} and \cite{Dkm,Mlf}.

\no {\bf Case $\hg_3\oplus\hg_3$ (type $(4,2)$)}.  The rational
forms of $\hg_3\oplus\hg_3$ are given by $\{\ngoq_k\}$, $k\geq 1$
square-free (see Proposition \ref{rath3h3}).  The fact that
$\ngoq_k$ is Anosov for any $k>1$ has been proved in several papers
(see \cite{Sml, Ito, AslSch, Mlf}) and it also follows from the
construction given in \cite{anosov} and Theorem \ref{gradedsum}. The
Pfaffian form of $\ngoq_1$ is $f_1(x,y)=x^2-y^2$, and thus it
follows from Proposition \ref{region}, (ii), that $\ngoq_1$ is not
Anosov.

\no {\bf Case $\ggo$ (type $(6,2)$)}.  It is proved in Section
\ref{cla} that the Lie algebra $\ggo$ defined in (\ref{62alg}) is
the only real Anosov Lie algebra of this type, and we have seen in
Proposition \ref{rat62} that $\ggo$ has only one rational form,
which is then the only rational Anosov Lie algebra of this type.

\no {\bf Case $\hg$ (type $(4,4)$)}.  We have seen in Section
\ref{cla} that the only possible real Anosov Lie algebras of this
type are the real forms of $\hg_{\CC}$, namely, $\hg$ and
$\ngoq_{-1}\otimes\RR$.  The rational forms of $\hg$ are determined
in Proposition \ref{rath}; they can be parametrized by $\hg^{\QQ}_k$
with $k$ a square-free natural number. We know that the Pfaffian
form of $\hg^{\QQ}_k$ is $f_k(x,y,z,w)=xw+y^2-kz^2$ and then
$Hf_k=4k$. By renaming the basis $\beta$ given in (\ref{base}) as
$\{ X_1,...,X_4,Z_1,...,Z_4\}$, we have that the Lie bracket of the
Anosov rational form $\hg^{\QQ}$ of $\hg$ defined by $\beta$ is
$$
\begin{array}{lcl}
[X_1,X_3]=Z_1+Z_3, && [X_2,X_3]=Z_2-Z_4, \\ \\

[X_1,X_4] = Z_2+Z_4, && [X_2,X_4]=(a^2-1)(Z_1-Z_3).
\end{array}
$$

This implies that the maps $J_Z$'s of $\hg^{\QQ}$ are given by
$$
J_{xZ_1+yZ_2+zZ_3+wZ_4}=\left[\begin{smallmatrix} 0&0&-x-z&-y-w\\ 0&0&-y+w&m(-x+z)\\ x+z&y-w&0&0\\
y+w&m(x-z)&0&0 \end{smallmatrix}\right],
$$
where $m=a^2-1$, and then its Pfaffian form is
$$
f(x,y,z,w)=mx^2-y^2-mz^2+w^2,
$$
with Hessian $Hf=16m^2$.  We know that $\hg^{\QQ}$ has to be
isomorphic to $\hg^{\QQ}_k$ for some square-free natural number $k$,
but in that case $f\simeq_{\QQ}f_k$ and so we would have
$Hf=q^2Hf_k$ for some $q\in\QQ^*$. Thus $16m^2=q^2k$, which implies
that $k=1$.  This shows that the Anosov rational forms of $\hg$
defined by different integers $a$'s are all isomorphic to
$\hg^{\QQ}_1$. In what follows, we shall prove that the other
rational forms of $\hg$ (i.e. $\hg^{\QQ}_k$ for $k>1$) are Anosov as
well.

Fix a square free natural number $k>1$.  Consider the basis
$\beta=\{ X_1,...,Z_4\}$ of $\hg^{\QQ}_k$ given in Proposition
\ref{rath} and set $\ngo_1=\la X_1,...,X_4\ra_{\QQ}$ and $\ngo_2=\la
Z_1,...,Z_4\ra_{\QQ}$.  Let $(a,b)\in\NN\times\NN$ any solution to
the Pell equation $x^2-ky^2=1$.  Let
$A:\hg^{\QQ}_k\mapsto\hg^{\QQ}_k$ be the linear map defined in terms
of $\beta$ by
\begin{equation}\label{auto44}
A_1=A|_{\ngo_1}=\left[\begin{smallmatrix} 0&0&b&-a\\ 0&0&-a&kb\\
0&1&2n&0\\ 1&0&0&2n \end{smallmatrix}\right], \qquad
A_2=A|_{\ngo_2}=\left[\begin{smallmatrix} 0&0&0&-1\\ 0&-a&b&4na\\
0&-bk&a&4nbk\\ -1&-2n&0&2n^2
\end{smallmatrix}\right].
\end{equation}
It is easy to check that $A\in\Aut(\hg^{\QQ}_k)$ for any $n\in\NN$,
and since $\det{A_1}=\det{A_2}=1$ we have that
$A_1,A_2\in\Sl_4(\ZZ)$, that is, $A$ is unimodular.  The
characteristic polynomial of $A_1$ is
$f(x)=(x^2-2nx+a-\sqrt{k}b)(x^2-2nx+a+\sqrt{k}b)$ and so its
eigenvalues are
\begin{equation}\label{eigen}
\begin{array}{l}
\lambda_1=n+\sqrt{n^2-a+\sqrt{k}b}, \qquad \lambda_2=n-\sqrt{n^2-a+\sqrt{k}b}, \\ \\
\mu_1=n+\sqrt{n^2-a-\sqrt{k}b}, \qquad
\mu_2=n-\sqrt{n^2-a-\sqrt{k}b}.
\end{array}
\end{equation}
We take $n\in\NN$ such that $a+\sqrt{k}b<n^2$.  Therefore
$1<\lambda_1$ and it follows from
$\lambda_1\lambda_2=a-\sqrt{k}b=\frac{1}{a+\sqrt{k}b}<1$ that
$\lambda_2<1$.  Also, $1<\mu_1$ and $\mu_1\mu_2=a+\sqrt{k}b>1$, and
hence $\mu_2\ne 1$, proving that $A_1$ is hyperbolic.  The
eigenvalues of $A_2$ are all of the form $\lambda_i\mu_j$. Indeed,
it can be checked that the eigenvector for $\lambda_i\mu_j$ is
$$
Z=Z_1-(a-\sqrt{k}b)\mu_jZ_2-(a+\sqrt{k}b)\lambda_iZ_3+\lambda_i\mu_jZ_4.
$$
Now, the fact that $\lambda_2<\mu_2<\mu_1<\lambda_1$ implies that
$\lambda_i\mu_j\ne 1$ for all $i,j$, showing that $A_2$ is also
hyperbolic and hence that $A$ is an Anosov automorphism of
$\hg^{\QQ}_k$.

The above is the most direct and shortest proof of the fact that
$\hg_k^{\QQ}$ is Anosov for any square free $k>1$, and it consists
in just checking that $A$ is unimodular and hyperbolic. But now we
would like to show where did this $A$ come from, which will show at
the same time that $\hg^{\QQ}_k$ is not Anosov for $k<0$.  Since the
proof of Proposition \ref{rath} actually shows that the set of
rational forms up to isomorphism of $\hg_{\CC}$ is given by
$$
\{\hg^{\QQ}_k:k\;\mbox{a nonzero square free integer number}\},
$$
this will prove that the real completion $\ngoq_{-1}\otimes\RR$ of
those with $k<0$ is not Anosov.

First of all, it is easy to see that any $\tilde{A}$ of the form
\begin{equation}\label{auto442}
\tilde{A}_1=\tilde{A}|_{\ngo_1}=\left[\begin{smallmatrix} B&0\\
0&C \end{smallmatrix}\right],\quad
\tilde{A}_2=\tilde{A}|_{\ngo_2}=\left[\begin{smallmatrix} b_{11}C&b_{12}C\\
b_{21}C&b_{22}C
\end{smallmatrix}\right]\quad \left(B=\left[\begin{smallmatrix} b_{11}&b_{12}\\ b_{21}&b_{22}
\end{smallmatrix}\right]\right),
\end{equation}
where $B,C\in\Gl_2(\CC)$, is an automorphism of $\hg_{\CC}$, for
which we are considering the basis $\alpha=\{X_1,...,Z_4\}$ with Lie
bracket defined as in (\ref{corh}).  Moreover, this forms a subgroup
of $\Aut(\hg_{\CC})$ containing the connected component at the
identity, since any other automorphism restricted to
$(\ngo_1)_{\CC}$ has the form $ \left[\begin{smallmatrix} 0&\star\\
\star&0
\end{smallmatrix}\right]$.  By taking $\tilde{A}^2$ if necessary, we can assume
that if $\hg^{\QQ}_k$ is Anosov then it admits an Anosov
automorphism of the form (\ref{auto442}).  The change of basis
matrix $P_k$ from the basis $\beta_k$ of the rational form
isomorphic to $\hg^{\QQ}_k$ given in the proof of Proposition
\ref{rath} to the basis $\alpha$ is
$$
P_k|_{\ngo_1}=\left[\begin{smallmatrix} \sqrt{k}&1&0&0\\
0&0&1&\sqrt{k}\\ -\sqrt{k}&1&0&0\\ 0&0&1&-\sqrt{k}
\end{smallmatrix}\right], \quad P_k|_{\ngo_2}=\left[\begin{smallmatrix} 2\sqrt{k}&0&0&0\\ 0&\sqrt{k}&-1&0\\
0&\sqrt{k}&1&0\\ 0&0&0&-2\sqrt{k} \end{smallmatrix}\right],
$$
and hence
$$
P_k^{-1}|_{\ngo_1}=\unm\left[\begin{smallmatrix} 1/\sqrt{k}&0&-1/\sqrt{k}&0\\ 1&0&1&0\\ 0&1&0&1\\
0&1/\sqrt{k}&0&-1/\sqrt{k} \end{smallmatrix}\right], \quad
P_k^{-1}|_{\ngo_2}=\unm\left[\begin{smallmatrix}
1/\sqrt{k}&0&0&0\\ 0&1/\sqrt{k}&1/\sqrt{k}&0\\ 0&-1&1&0\\
0&0&0&-1/\sqrt{k} \end{smallmatrix}\right].
$$
We then have that $A=P_k^{-1}\tilde{A}P_k\in\Aut(\hg^{\QQ}_k)$ if
and only if $A_1=P_k^{-1}\tilde{A}_1P_k$ and
$A_2=P_k^{-1}\tilde{A}_2P_k$ belong to $\Gl_4(\QQ)$.  A
straightforward computation shows that $A_1,A_2\in\Gl_4(\ZZ)$ (i.e.
$A$ is unimodular) if and only if $B\in\Gl_2(\ZZ[\sqrt{k}])$,
$C=\overline{B}$ and $\det{B}\overline{\det{B}}=\pm 1.$  Here,
$\ZZ[\sqrt{k}]$ is the integer ring of the quadratic numberfield
$\QQ[\sqrt{k}]$ and the conjugation is defined, as usual, by
$\overline{x+\sqrt{k}b}=x-\sqrt{k}b$ for all $x,y\in\QQ$. Recall
that if $\det{B}=a-\sqrt{k}b$, $a,b\in\ZZ$, and we assume that
$\det{B}\overline{\det{B}}=1$, then $a^2-kb^2=1$, the Pell equation.
In order to make easier the computation of eigenvalues we can take
$B$ in its rational form, say
$$
B=\left[\begin{smallmatrix} 0&-a+\sqrt{k}b\\ 1&2n
\end{smallmatrix}\right], \qquad
\overline{B}=\left[\begin{smallmatrix} 0&-a-\sqrt{k}b\\ 1&2n
\end{smallmatrix}\right],
$$
for some $n\in\ZZ$.  This implies that the characteristic polynomial
of $\tilde{A}_1$ is
$f(x)=(x^2-2nx+a-\sqrt{k}b)(x^2-2nx+a+\sqrt{k}b)$ and so the
eigenvalues of $\tilde{A}_1$ and $A_1$ are as in (\ref{eigen}).
Concerning the hyperbolicity, if $k<0$ then either $b=0$ or $a=0$
and $k=-1$, which in any case implies that $|\mu_1\mu_2|=1$, a
contradiction.  Therefore $\hg^{\QQ}_k$ is not Anosov for $k<0$, as
was to be shown.  For $k>0$, we can easily see that conditions
$a,b,n\in\NN$, $a+\sqrt{k}b<n^2$, are enough for the hyperbolicity
of $A_1$. For $A_2$, we can use the following general fact: the
eigenvalues of a matrix of the form $\tilde{A}_2$ in (\ref{auto442})
are precisely the possible products between eigenvalues of $B$ and
eigenvalues of $C$; and so the hyperbolicity of $A_2$ follows as in
the short proof.

We finally note that $A=P_k^{-1}\tilde{A}P_k$ with this $B$ is
precisely the automorphism proposed in (\ref{auto44}).

{\small
\begin{table}
$$
\begin{array}{cccccc}
\hline\hline
&&&&& \\
Real \, Anosov & Dimension & Type & Anosov  & Non-Anosov & Signature  \\
Lie \, algebra & &  & rat. \, forms & rat. \, forms & \\ \\
\hline
&&&&&\\
\RR^n,\,2\le n\le8 & n & n & \QQ^n & - - & {\rm any} \\ \\

\hg_3\oplus\hg_3   &6 &(4,2) & \ngoq_k,\;k>1&\ngoq_1 & \{ 3,3\} \\ \\

\fg_3             &6 &(3,3) &\fg_3^\QQ&  - - &  \{ 3,3\} \\ \\

\hg_3\oplus\hg_3\oplus\RR^2 & 8 &(6,2) & \ngoq_k\oplus\QQ^2,\;k>1&\ngoq_1\oplus\QQ^2 & \{ 4,4\} \\ \\

\fg_3\oplus\RR^2 & 8 &(5,3) &\fg_3^\QQ\oplus \QQ^2 & - - & \{ 4,4\}  \\ \\

\ggo& 8 &(6,2)& \ggo^\QQ& - - & \{ 4,4\} \\ \\

\hg & 8& (4,4) & \hg_k^\QQ,\; k\ge 1& - - & \{ 4,4\}  \\ \\

\lgo_4\oplus\lgo_4 &8& (4,2,2) & \lgo_k^\QQ,\; k>1& \lgo_1^\QQ & \{ 4,4\} \\ \\

\hline \hline \\
\end{array}
$$
\caption{Real and rational Anosov Lie algebras of dimension $\le
8$.}\label{Anosovtable}
\end{table}}

\begin{remark}
{\rm An alternative proof of the fact that any rational form of
$\hg$ is Anosov can be given by using \cite[Corollary 2.3]{Dn}.
Indeed, the subgroup
$$
S=\Sl_2(\RR)\times\Sl_2(\RR)=\left\{A\in\Aut(h):A_1=\left[\begin{smallmatrix} B&0\\
0&C \end{smallmatrix}\right], \quad B,C\in\Sl_2(\RR)\right\}
$$
is connected, semisimple and all its weights on $\hg$ are
non-trivial.  Recall that such a corollary can not be applied to the
cases $\hg_3\oplus\hg_3$ and $\lgo_4\oplus\lgo_4$, as they admit a
rational form which is not Anosov.}
\end{remark}

\no {\bf Case $\lgo_4\oplus\lgo_4$ (type $(4,2,2)$)}.  The rational
forms of $\lgo_4\oplus\lgo_4$ are determined in Proposition
\ref{ratl4l4} and they are denoted by $\lgo_k^\QQ$, $k$ a square
free natural number.  Let $\beta$ denote the basis of $\lgo_k^\QQ$
given in the proposition. For $a \in \ZZ,$ $a \ge 2,$ consider the
hyperbolic matrix
$$
B=\left[\begin{smallmatrix}a & a^2-1\\1& a
\end{smallmatrix}\right] \in \Sl_2(\ZZ),
$$
with eigenvalues $\l_1=a+(a^2-1)^\unm$ and $\l_2=a-(a^2-1)^\unm.$ It
is easy to check that the linear map $A:\lgo_b^\QQ \longrightarrow
\lgo_b^\QQ$ whose matrix in terms of $\beta$ is
$$
[A]_\beta=\left[\begin{array}{ccc}B &&\\ &
\ddots&\\&&B\end{array}\right]
$$
is an automorphism of $\lgo_b^\QQ$ for $b=a^2-1.$  $A$ is hyperbolic
since $\l_1>1>\l_2$ and it is unimodular by definition, so that $A$
is an Anosov automorphism.  Recall that $\lgo_k^\QQ \simeq
\lgo_{k'}^\QQ$ if and only if $k=q^2k'$ for some $q \in \QQ^*$ (see
Proposition \ref{ratl4l4}).  Given a square-free natural number $k
>1,$ there always exist $a,q \in \ZZ$ such that  $a^2-1=q^2k$
(Pell equation), and thus any $\lgo_k^\QQ$ with $k>1$ square free is
Anosov.

We now prove that $\lgo_1^\QQ$ is not Anosov. In the proof of
Proposition \ref{ratl4l4} we have showed that any $A \in
\Aut(\lgo_1^\QQ)$ has the form (\ref{Aform}) and satisfies
$$
q f(z,w) = f(A^t_4(z,w)) \qquad \qquad \forall (z,w) \in \QQ^2,
$$
where $q=\det A_3A_1$ and $f(z,w) = z^2-w^2.$ In the same spirit of
Proposition \ref{region}, this implies that $A_4^t$ leaves a finite
set invariant and so it can never be hyperbolic.

The results obtained in this section can be summarized as follows.

\begin{theorem}\label{clasrat}
Up to isomorphism, the rational Anosov Lie algebras of dimension
$\le 8$ are

{\rm \begin{itemize} \item $\QQ^n,$ \,$n=2,\dots,8,$ \; ($\RR^n$),

\item $\ngoq_k,$ \,$k \ge 2,$ \;  ($\hg_3\oplus\hg_3$),

\item $\fg_3^\QQ,$ \;  ($\fg_3$),

\item $\ngoq_k\oplus \QQ^2,$ \,$k \ge 2,$ \; ($\hg_3\oplus\hg_3\oplus \RR^2$),

\item $\fg_3^\QQ \oplus\QQ^2,$ \; ($\fg_3\oplus\RR^2$),

\item $\ggo^\QQ,$ \; ($\ggo$),

\item $\hg_k^\QQ,$ \,$k \ge 1,$ \; ($\hg$),

\item $\lgo_k^\QQ,$ \,$k \ge 2,$ \;($\lgo_4\oplus\lgo_4$),

\end{itemize}}

\no where $k$ always run over square-free numbers and the Lie
algebra between parenthesis is the corresponding real completion.
\end{theorem}

In the last column of Table \ref{Anosovtable} appear the signatures
of the Anosov automorphisms found in each case.  It follows from the
proofs given in Section \ref{cla} that the eigenvalues of any Anosov
automorphism always appear in pairs $\{\l,\l^{-1}\}$ (with only one
exception: $\fg_3$), and thus there is only one possible signature
for each nonabelian Anosov Lie algebra of dimension $\leq 8$.

\begin{corollary}\label{sgn}
Let $N/\Gamma$ be a nilmanifold (or infranilmanifold) of dimension
$\leq 8$ which admits an Anosov diffeomorphism.  Then $N/\Gamma$ is
either a torus (or a compact flat manifold) or the dimension is $6$
or $8$ and the signature is $\{ 3,3\}$ or $\{ 4,4\}$, respectively.
\end{corollary}

It is not true in general that there is only one possible signature
for a given Anosov Lie algebra.  For instance, it is easy to see
that the free $2$-step nilpotent Lie algebra on $4$ generators
admits Anosov automorphisms of signature $\{ 4,6\}$ and $\{ 5,5\}$.

\section{Appendix: Algebraic numbers}\label{algnum}

We will give in the following a short summary of some results about
algebraic numbers over $\QQ$ that are used throughout the
classification. We are mainly following \cite[Chapter V]{Lng}. Note
that we will omit information on numberfields since we are not going
to need it.

An element $\l \in \CC$ is called {\it algebraic over $\QQ$} if
there exist a polynomial $p(x) \in \QQ[x]$ such that $p(\l)=0$. It
is easy to see that the set $D$ of all such polynomials  form an
ideal in $\QQ[x]$ and since this is a principal ideal domain, $D$ is
generated by a single polynomial. This polynomial
 can be chosen to be monic, and in that case it is uniquely determined by $\l$ and
 will be called the {\it minimal
 polynomial of $\l$}, denoted by $m_\l(x)$.  Therefore, if we have an algebraic number $\l$
 then we can define {\it the degree
 of $\l$} as the degree of $m_\l(x)$.  It will be denoted by $\degr\l$.  The minimal
 polynomial $m_\l(x)$ is  irreducible over $\QQ$ and $\l$
 is not a double root of $m_\l(x)$.

If $\l\ne\mu$ are two algebraic numbers, we say that they are {\it
conjugated} if $m_\l(\mu)=0$.  Note that the numbers which are
conjugated to $\l$ are uniquely determined by $\l$ and have the same
degree.

An algebraic number $\l$ is said to be an {\it algebraic integer} if
there exists a monic polynomial $p(x)\in\ZZ[x]$ such that $p(\l)=0$.
It can be seen that in this case,  $m_\l(x) \in \ZZ[x]$, and
moreover, these conditions are actually equivalent. An algebraic
number is called {\it totally real} if $m_\l(x)$ has only real
roots, that is, $m_\l(x) = \displaystyle{\prod_{i=1}^r} (x-\l_i)$
with   $\l_i \in \RR, \; \l_1 =\l.$ If $\l$ is a totally real
algebraic number with $\degr \l =r$, set
$A_\l = \left[\begin{smallmatrix} \l_1 & & \\
& \ddots& \\&& \l_r \end{smallmatrix}\right]$.  The characteristic
polynomial of $A_\l$ is $m_\l(x)$ and then the rational form of
$A_\l$ is given by
$$
\left[\begin{smallmatrix} 0&0 & \dots &0 &-a_0 \\
1& 0 & \dots& 0& -a_1 \\ 0&1& &0 & -a_2\\ \vdots & & \ddots &
&\vdots  \\ 0& & & 1& -a_{r-1}
\end{smallmatrix}\right],
$$
where $m_\l(x)= x^r+a_{r-1}x^{r-1}+ \dots + a_1x+a_0$.  If $\l$ is
an algebraic integer then $a_i \in \ZZ $ for all $i =0, \dots, r-1$
and then this shows that $A_\l$ is conjugated to an element in
$\Gl_r(\ZZ)$.

Conversely, if $A = \left[\begin{smallmatrix} \l_1 & & \\ & \ddots&
\\&& \l_r
\end{smallmatrix}\right]$ is conjugated to an element of $\Gl_r(\ZZ)$, then if
$p_A(x)$ is the characteristic polynomial of $A$, $p_A(x) \in
\ZZ[x]$, and therefore $\l_i$ is an algebraic integer for all $i=1,
\dots,r$.  Concerning the degree of the $\l_i$'s as algebraic
numbers in such a case, we can only say that $1 \le \degr \l_i \le
r.$ Moreover, if $\l_i=\l_j$ for some $i \ne j$, and since $\l$ is
not a double root of $m_\l(x)$, we will have that
$m_{\l_i}^2(x)|p_A(x)$ and hence $1\le 2 \degr \l_i \le r.$

If $\l$ is an algebraic integer, we say that $\l$ is a {\it unit} if
$1/\l$ is an algebraic integer as well. If it is so, then the
constant coefficient $a_0$ of $m_\l(x)$ is $(- 1)^n$, where $n=\degr
\l$. Conversely, if $a_0=\pm 1$ then $\l$ is a unit.


\begin{thebibliography}{99}

\bibitem{AslSch} {\sc L. Auslander, J. Scheuneman}, On certain automorphisms of nilpotent
Lie groups, {\it Global Analysis: Proc. Symp. Pure Math.} {\bf 14}
(1970), 9-15.

\bibitem{BnsLbr} {\sc Y. Benoist, F. Labourie},  Sur le diffeomorphismes d'Anosov affines a feuilletages stable et instable differentiables, {\it Inventiones Math.} {\bf 111} (1993), 285-308.

\bibitem{CssKnnScv} {\sc C. Cassidy, N. Kennedy, D. Scevenels}, Hyperbolic automorphisms for groups in $\tca(4,2)$,
{\it Contemporary Math.} {\bf 262} (2000), 171-175.

\bibitem{Dn} {\sc S.G. Dani}, Nilmanifolds with Anosov automorphism, {\it J. London
Math. Soc.} {\bf 18} (1978), 553-559.

\bibitem{Dkm} {\sc K. Dekimpe}, Hyperbolic automorphisms and Anosov diffeomorphisms on
nilmanifolds, {\it Trans. Amer. Math. Soc.} {\bf 353} (2001),
2859-2877.

\bibitem{DkmDsc} {\sc K. Dekimpe, S. Deschamps}, Anosov diffeomorphisms on a class of
$2$-step nilmanifolds, {\it Glasg. Math. J.} {\bf 45} (2003),
269-280.

\bibitem{DkmMlf} {\sc K. Dekimpe, W. Malfait}, A special class of nilmanifolds admitting
an Anosov diffeomorphism, {\it Proc. Amer. Math. Soc.} {\bf 128}
(2000), 2171-2179.

\bibitem{Ebr} {\sc P. Eberlein}, Geometry of $2$-step nilpotent Lie groups, preprint
2003 (author's web page).

\bibitem{Frd} {\sc D. Fried}, Nontoral pinched Anosov maps, {\it Proc. Amer. Math. Soc.}
{\bf 82} (1981), 462-464.

\bibitem{Frn} {\sc J. Franks}, Anosov diffeomorphisms, {\it Global Analysis:
Proc. Symp. Pure Math.} {\bf 14} (1970), 61-93.

\bibitem{Ggr} {\sc M. Gauger}, On the classification of metabelian Lie algebras,
{\it Trans. Amer. Math. Soc.} {\bf 179} (1973), 293-329.

\bibitem{Grm} {\sc M. Gromov}, Groups of polynomial growth and expanding maps,
{\it Inst. des Hautes Etudes Sci.} {\bf 53} (1981), 53-73.

\bibitem{GrnSgl1} {\sc F. Grunewald, D. Segal}, Nilpotent groups of Hirsch lenght six,
{\it Math. Z. } {\bf 179} (1982), 162-175.

\bibitem{GrnSgl2} {\sc F. Grunewald, D. Segal}, Reflections on the classification of
torsion-free nilpotent groups, {\it Group Theory: essays for Philip
Hall} (1984), 159-206, Academic Press.

\bibitem{Ito} {\sc K. Ito}, Classification of nilmanifolds $M^n$ $(n\leq 6)$
admitting Anosov diffeomorphisms, {\it The study of dynamical
systems} Kyoto (1989), 31-49, World Sci. Adv. Ser. Dynam. Systems
{\bf 7}.

\bibitem{KS} {\sc A. Katok, R. Spatzier}, Differential rigidity of Anosov actions of higher
rank abelian groups and algebraic lattice actions, {\it Proc.
Steklov Inst. Math} {\bf 216} (1997), 287-314.

\bibitem{Lng} {\sc S. Lang}, Algebra, {\it Addison-Wesley} 1993.

\bibitem{anosov} {\sc J. Lauret}, Examples of Anosov diffeomorphisms, {\it Journal of Algebra} {\bf 262} (2003), 201-209.  {\it Corrigendum}: {\bf 268} (2003), 371-372.

\bibitem{Mgn} {\sc L. Magnin}, Sur les algebres de Lie nilpotents de dimension $\leq$ 7,
{\it J. Geom. Phys.}, 1986 Vol. {\bf III}, 119-144.

\bibitem{Mlf} {\sc W. Malfait}, Anosov diffeomorphisms on nilmanifolds of
dimension at most $6$, {\it Geom. Dedicata} {\bf 79} (2000),
291-298.

\bibitem{Mlf2} {\sc W. Malfait}, An obstruction to the existence of Anosov diffeomorphisms
on infra-nilmanifolds, {\it Contemp. Math.} {\bf 262} (2000),
233-251 (Kortrijk 1999).


\bibitem{Mnn} {\sc A. Manning}, There are no new Anosov  diffeomorphisms on tori,
{\it Amer. J. Math.} {\bf 96} (1974), 422-429.

\bibitem{Mrg} {\sc G. Margulis}, Problems and conjectures in rigidity theory, {\it Mathematics: Frontiers and perspectives 2000}, IMU.

\bibitem{Prt} {\sc H. Porteous}, Anosov diffeomorphisms of flat manifolds,
{\it Topology} {\bf 11} (1972), 307-315.

\bibitem{Rgn} {\sc M.S. Raghunathan}, Discrete subgroups of Lie groups,
{\it Ergeb. Math.} {\bf 68} (1972), Springer-Verlag.

\bibitem{Sch} {\sc J. Scheuneman}, Two-step nilpotent Lie algebras, {\it J. Algebra} {\bf 7} (1967), 152-159.

\bibitem{Sly} {\sc C. Seeley}, $7$-dimensional nilpotent Lie algebras, {\it Trans. Amer. Math. Soc.},
{\bf 335} (1993), 479-496.

\bibitem{Shb} {\sc M. Shub}, Endomorphisms of compact differentiable manifolds,
{\it Amer. J. Math.} {\bf 91} (1969), 175-199.

\bibitem{Sml} {\sc S. Smale}, Differentiable dynamical systems,
{\it Bull. Amer. Math. Soc.} {\bf 73} (1967), 747-817.

\bibitem{Smn} {\sc Y. Semenov}, On the rational forms of nilpotent Lie algebras and lattices in nilpotent Lie groups,
{\it L'Enseignement Mathematique} {\bf 48} (2002), 191-207.

\bibitem{Vrj} {\sc A. Verjovsky}, Sistemas de Anosov, Notas de Curso, IMPA (Rio de Janeiro).

\end{thebibliography}
\end{document}